\documentclass[11pt,pdftex]{amsart}
\usepackage[a4paper,margin=30mm]{geometry}
\usepackage{amsmath,amssymb,amsthm,amscd,mathtools}
\usepackage{mathrsfs}
\usepackage{tikz-cd}
\usepackage[colorlinks=true,linkcolor=blue,citecolor=blue,urlcolor=blue]{hyperref}
\usepackage[msc-links,backrefs]{amsrefs}

\numberwithin{equation}{section}

\theoremstyle{plain}
\newtheorem{thm}{Theorem}[section]

\newtheorem{pro}[thm]{Proposition}
\newtheorem{cor}[thm]{Corollary}

\theoremstyle{definition}
\newtheorem{defn}[thm]{Definition}
\newtheorem{eg}[thm]{Example}

\theoremstyle{remark}
\newtheorem{rem}[thm]{Remark}

\newcommand{\Image}{{\operatorname{Im}}}

\newcommand{\Ker}{{\operatorname{Ker}}}
\newcommand{\Kernel}{{\operatorname{Ker}}}
\newcommand{\Conj}{\operatorname{Conj}}
\newcommand{\As}{\operatorname{As}}
\newcommand{\Inn}{\operatorname{Inn}}
\newcommand{\Aut}{\operatorname{Aut}}
\newcommand{\Stab}{\operatorname{Stab}}
\newcommand{\Alex}{\operatorname{Alex}}

\newcommand{\XMod}{\mathsf{XMod}}
\newcommand{\AugQdl}{\mathsf{AugQdl}}
\newcommand{\AugRack}{\mathsf{AugRack}}

\begin{document}

\title[Freeness of Crossed Modules and Augmented Quandles]
{Freeness of crossed modules, augmented quandles, \\
and augmented racks}
\author{Toshiyuki Akita}

\address{Department of Mathematics, Faculty of Science, Hokkaido University,
Sapporo, 060-0810 Japan}
\email{akita@math.sci.hokudai.ac.jp}

\author{Sota Takase}

\address{Department of Mathematics, Graduate School of Science, Hokkaido University,
Sapporo, 060-0810 Japan}

\curraddr{\textsc{SMBC Global Investment \& Consulting Ltd., 
Tokyo, 105-6427 Japan}}
\email{takase.sota@smbc-gic.co.jp}


\date{}
\begin{abstract}
We study freeness of crossed modules, augmented
quandles, and augmented racks. 
We show that the associated quandle functor and the associated crossed module functor 
are left adjoint to the natural
forgetful functors.
For a map $f\colon S\to G$, we define the free augmented rack
$R_f\to G$ and the free augmented quandle $Q_f\to G$. 
If $C_f\to G$ denotes the free crossed module determined by $f$,
then the associated crossed modules of $R_f\to G$ and $Q_f\to G$
are both isomorphic to $C_f\to G$.
We also characterize free augmented quandles in terms of stabilizers.
As applications, we prove that the canonical inner
augmentation of every generalized Alexander quandle
is free. For a classical knot quandle, the canonical
inner augmentation is free if and only if the knot is
a torus knot. For every $n$-knot with $n\geq2$, the
canonical inner augmentation is also free.
\end{abstract}
\maketitle

\section{Introduction}

J.~H.~C. Whitehead introduced crossed modules and free crossed modules
in his study of low-dimensional homotopy groups
\cites{MR0004123,MR0017537,MR0030760}.
Crossed modules provide algebraic models for connected homotopy $2$-types
~\cite{MR2841564}, and crossed modules are equivalent to categorical
groups (strict $2$-groups in modern terminology)~\cite{MR0419643}.
Crossed modules also appear in higher gauge theory
~\cites{MR2825807,MR2342821} and higher holonomy~\cite{Kohno2020}, and
crossed modules of Lie algebras have been studied in Lie theory.
See \cites{MR4322147,MR4368859} for example.

Free crossed modules arise naturally from relative homotopy groups.
For a pair of path-connected spaces $(X,A)$, the boundary map
\[
\partial\colon\pi_2(X,A)\to\pi_1(A)
\]
is a crossed module, called the fundamental crossed module of the pair $(X,A)$.
Whitehead proved that if $X$ is obtained from $A$ by attaching
$2$-cells, then the fundamental crossed module is free on the attached
$2$-cells~\cite{MR0030760}*{Section~16}. Thus free crossed modules arise
naturally in topology and in the theory of identities among relations
for group presentations~\cite{BrownHuebschmann1982}.

Quandles were introduced independently by Joyce and Matveev in connection with knot theory.
See \cites{MR2628474,MR0638121,MR0672410}.
Racks generalize quandles. 
The notion of a rack predates
that of a quandle, but the term ``rack'' was introduced
by Fenn and Rourke, who developed the theory systematically
and applied it to knots and links in codimension two~\cite{MR1194995}.
Racks and quandles also yield set-theoretic solutions of the
Yang--Baxter equation~\cites{MR1994219,arXiv:2506.23175}. 
Moreover,
racks and quandles equipped with suitable
$2$-cocycles give rise to braided vector spaces and Nichols algebras,
which play an important role in the study of pointed Hopf algebras
~\cite{MR1994219}.

Augmented racks and augmented quandles are racks and quandles equipped
with compatible group actions and augmentation maps. A crossed module is an augmented quandle,
and an augmented quandle is an augmented rack.
We write $\XMod$, $\AugQdl$, and $\AugRack$ for the corresponding
categories. Let
\[
U_{\mathrm Q}\colon\XMod\to\AugQdl,
\quad
U_{\mathrm R}\colon\AugQdl\to\AugRack
\]
denote the forgetful functors.

The associated quandle construction for augmented racks is
standard, while Fenn and Rourke constructed the associated
crossed module of an augmented rack~\cite{MR1194995}*{pp.~355--356}.
In this paper, we show that the associated quandle functor and
the associated crossed module functor are left adjoint to the
corresponding forgetful functors. 
Thus we obtain the adjunctions
\[
\AugRack
\underset{U_{\mathrm R}}{\stackrel{Q}{\rightleftarrows}}
\AugQdl
\underset{U_{\mathrm Q}}{\stackrel{\As}{\rightleftarrows}}
\XMod,
\qquad
Q\dashv U_{\mathrm R},
\quad
\As\dashv U_{\mathrm Q}.
\]
For all three categories, the free objects considered in this paper are determined by a map
$f\colon S\to G$, where $S$ is a set indexing the basis elements and $f$ records their images in $G$.
Let $R_f\to G$, $Q_f\to G$, and $C_f\to G$ denote, respectively, the free augmented rack, the free augmented quandle, and the free crossed module determined by the map $f$.
We give explicit
descriptions of $R_f$ and $Q_f$ and compare their associated crossed modules
with $C_f\to G$.

Theorem~\ref{thm:AT3} shows that the associated crossed modules of
$R_f\to G$ and $Q_f\to G$ are both isomorphic to $C_f\to G$ via isomorphisms 
preserving the corresponding bases.
Thus the two left adjoints above preserve the corresponding free objects
determined by the same map $f$. 
In particular, every free crossed module is the associated crossed
module of both a free augmented rack and a free augmented quandle.

The associated crossed module of a free augmented quandle is free, but
the converse does not hold in general.
 We describe the additional condition required for an augmented
quandle to be free in terms of stabilizers. 
Let $\partial\colon Q\to G$ be an augmented quandle. 
Choose a complete set $S\subseteq Q$ of representatives for the orbits of the $G$-action on $Q$, let
$i\colon S\hookrightarrow Q$ be the inclusion, and put
$f=\partial\circ i$. By Proposition~\ref{pro:free-aq-coset-criterion},
$\partial$ is free on $i$ if and only if
\[
\Stab_G(i(s))=\langle f(s)\rangle
\quad(s\in S).
\]
Since $\langle f(s)\rangle\leq\Stab_G(i(s))$ always holds,
freeness is characterized by equality.

For a generalized Alexander quandle $M$, we show that the canonical inner augmentation
$M\to\Inn(M)$ is a free augmented quandle, with a basis given by a complete set of
representatives for the $\Inn(M)$-orbits
(Theorem~\ref{thm:generalized-Alexander-free}).
Knot quandles show that freeness of an augmented quandle and freeness
of its associated crossed module need not coincide.
 If $K$ is a nontrivial classical knot,
then the augmentation $Q_K\to G_K$ is not a free augmented quandle,
although its associated crossed module
is free. 
On the other hand, the canonical inner augmentation $Q_K\to\Inn(Q_K)$ is free if and only if
$K$ is a torus knot
(Theorem~\ref{thm:one-knot-free-augmented-quandle}).
In contrast, for every $n$-knot with $n\geq2$, the canonical inner
augmentation is free
(Theorem~\ref{thm:higher-knot-free-augmented-quandle}).
Thus freeness of the associated crossed module
is necessary but not sufficient for freeness of
the augmented quandle.

The paper is organized as follows.
Section~\ref{sec:crossed-modules} reviews crossed modules and free
crossed modules, including Whitehead's theorem and Ratcliffe's
criterion.
Section~\ref{sec:augmented-quandles-racks} recalls racks, quandles,
augmented racks, augmented quandles, and the associated quandle
construction.
Section~\ref{sec:associated-crossed-modules} describes the associated
crossed module construction and establishes the adjunctions above.
In Section~\ref{sec:free-augmented}, we introduce free augmented racks
and free augmented quandles and compare their associated crossed modules
with free crossed modules determined by the same map $f\colon S\to G$.
Section~\ref{sec:stabilizers} gives the stabilizer description and the
freeness criterion.
Section~\ref{sec:examples-applications} gives examples and applications,
including conjugation quandles, generalized Alexander quandles, and knot
quandles.

Some of the results were announced in~\cite{AkitaTakaseOCAMI2025}.
Part of this work is based on the second author's master's
thesis~\cite{Takase2024}, written under the supervision of the first
author.

\subsection*{Notation and conventions}

All group actions in this paper are right actions.
All $G$-modules are right $\mathbb ZG$-modules.
For a set $X$, we denote by $F(X)$ the free group on $X$.
For a group $G$ and an element $g\in G$, we write
$\langle g\rangle$ for the cyclic subgroup generated by $g$.
For a subset $A\subset G$, we denote the
centralizer of $A$ in $G$ by $C_G(A)$.
For $g\in G$, we write $C_G(g)$ for $C_G(\{g\})$.
The center of $G$ is denoted by
$Z(G)$.
If a group $G$ acts on a set $X$, the stabilizer (or isotropy
subgroup) of $x\in X$ is denoted by
\[
\Stab_G(x)\coloneqq
\{g\in G\mid x\cdot g=x\}.
\]

\section{Crossed modules and free crossed modules}\label{sec:crossed-modules}
\subsection{Definitions and examples}
Crossed modules were introduced by J.~H.~C.~Whitehead in his study of low-dimensional homotopy groups
\cites{MR0004123,MR0017537,MR0030760}.

\begin{defn}
A \emph{crossed module} consists of a group homomorphism
$\partial\colon C\to G$
and a right action of $G$ on $C$ by group automorphisms such that
\begin{enumerate}
\item
$\partial(c\cdot g)=g^{-1}(\partial c)g$
\quad$(c\in C,g\in G)$,
\item
$c_{1}\cdot\partial(c_{2})=c_{2}^{-1}c_{1}c_{2}$
\quad$(c_{1},c_{2}\in C)$.
\end{enumerate}
The homomorphism $\partial$ is called the \emph{boundary map}. The identities in (2) are called the Peiffer identities. A crossed module with codomain $G$ will be called a \emph{crossed $G$-module}.
\end{defn}

If $\partial\colon C\to G$ is a crossed module, then the following properties hold.
\begin{enumerate}
\item
$N\coloneqq\Image\,\partial$ is a normal subgroup of $G$.
\item
$A\coloneqq\Kernel\,\partial$ is a central subgroup of $C$.
\item
The $G$-action on $C$ induces a $G/N$-action on $A$.
\end{enumerate}

\begin{defn}
A \emph{morphism of crossed modules}
$(\beta,\alpha)\colon (\partial \colon C \to G) \to (\delta \colon D \to H)$
is a pair of group homomorphisms
$\alpha \colon G \to H$ and $\beta \colon C \to D$
such that the following diagram commutes:
\begin{equation}\label{eq:morphism}
\begin{tikzcd}
C \times G \arrow[r, ""] \arrow[d, "\beta \times \alpha"'] &
C \arrow[r, "\partial"] \arrow[d, "\beta"'] &
G \arrow[d, "\alpha"'] \\
D \times H \arrow[r, ""] &
D \arrow[r, "\delta"] &
H
\end{tikzcd}
\end{equation}
The horizontal maps on the left are the corresponding group actions. A morphism $(\beta,\alpha)$ is an \emph{isomorphism} if both $\beta$ and $\alpha$ are group isomorphisms. In this case $(\beta^{-1},\alpha^{-1})$ is also a morphism of crossed modules and is the inverse of $(\beta,\alpha)$.
\end{defn}

We recall some basic examples of crossed modules.
\begin{eg}[J.~H.~C.~Whitehead~\cite{MR0030760}]
For a pair $(X,A)$ of path-connected spaces, the boundary map
\[
\partial\colon \pi_2(X,A)\to \pi_1(A)
\]
is a crossed module. It is called the \emph{fundamental crossed module} of the pair $(X,A)$. Conversely, every crossed module can be realized as the fundamental crossed module of a CW-pair. See Brown--Higgins--Sivera~\cite{MR2841564}*{Proposition 2.4.7} for this realization and for the description of homotopy $2$-types by crossed modules.
\end{eg}

\begin{eg}[Quillen]
For a Serre fibration $F\to E\to B$ of path-connected spaces, the induced homomorphism
$\pi_1(F)\to\pi_1(E)$
is a crossed module.
See Brown--Higgins--Sivera~\cite{MR2841564}*{Section~2.6}.
\end{eg}

\begin{eg}
If $N \lhd G$ is a normal subgroup, then the inclusion homomorphism
$N \hookrightarrow G$
is a crossed module.
\end{eg}

\begin{eg}
A central extension
$\partial\colon E \twoheadrightarrow G$
is a crossed module with the right $G$-action on $E$ induced by conjugation in $E$.
\end{eg}

\begin{eg}
If $M$ is a $G$-module, then the zero homomorphism
$M \to G$
is a crossed module.
\end{eg}

\subsection{Free crossed modules}
Whitehead introduced free crossed modules in connection with Theorem~\ref{thm:Whitehead} below~\cite{MR0030760}.

\begin{defn}
Let $\partial \colon C \to G$ be a crossed module, let $S$ be a set, and let $\phi \colon S \to C$ be a map.
We say that $\partial \colon C \to G$ is \emph{free on $\phi\colon S\to C$} if the following universal property holds. For every crossed module
$\delta \colon D \to H$, every group homomorphism
$\alpha \colon G \to H$, and every map $\psi \colon S \to D$ satisfying
$\alpha \circ \partial \circ \phi = \delta \circ \psi$, there exists a unique group homomorphism
$\beta \colon C \to D$
such that $\beta \circ \phi = \psi$ and the following diagram is a morphism of crossed modules.
\begin{equation}\label{eq:free}
\begin{tikzcd}
S \arrow[r, "\phi"] \arrow[rd, "\psi"'] &
C \arrow[r, "\partial"] \arrow[d, "\exists!\,\beta", dotted] &
G \arrow[d, "\alpha"] \\
& D \arrow[r, "\delta"] & H .
\end{tikzcd}
\end{equation}
\end{defn}

In this situation, $\phi\colon S\to C$ is injective. Indeed, let $\mathbb Z[S]$ be the free abelian group on $S$ and consider the zero crossed module $\mathbb Z[S]\to 1$. Apply the universal property to the unique homomorphism $G\to1$ and to the basis inclusion $\iota\colon S\hookrightarrow\mathbb Z[S]$. We obtain a homomorphism $\beta\colon C\to\mathbb Z[S]$ such that $\beta\circ\phi=\iota$. Since $\iota$ is injective, so is $\phi$. We may therefore identify $S$ with its image in $C$. In this case, $S$ is called a \emph{basis} of the free crossed module $\partial\colon C\to G$.

The definition above is equivalent to the usual definition of a
free crossed $G$-module on the map
$\partial\circ\phi\colon S\to G$ in the category
$\XMod_G$, whose morphisms are over $\mathrm{id}_G$.
One direction is immediate by taking $\alpha=\mathrm{id}_G$. 
Conversely, suppose that the universal property holds for crossed $G$-modules.
Given a crossed module $\delta\colon D\to H$ and a group homomorphism
$\alpha\colon G\to H$, pull back $\delta$ along $\alpha$ to a crossed
$G$-module and apply the universal property there.
This gives the unique morphism required in the definition above.

The definition involves the basis map $\phi\colon S\to C$, but the isomorphism type of a free crossed module is determined by the composite
\[
f=\partial\circ\phi\colon S\to G.
\]
More precisely, we have the following proposition.

\begin{pro}\label{pro:free-CM-determined-by-f}
Let $G$ be a group, let $S$ be a set, and let $f\colon S\to G$ be a map. For $i=1,2$, suppose that $\partial_i\colon C_i\to G$ is free on a map $\phi_i\colon S\to C_i$ and that $\partial_i\circ\phi_i=f$. Then there exists a unique isomorphism $\beta\colon C_1\to C_2$ making the following diagram commute.
\[
\begin{tikzcd}[column sep=large]
S \arrow[r,"\phi_1"] \arrow[dr,"\phi_2"'] &
C_1 \arrow[d,"\beta","\cong"'] \arrow[r,"\partial_1"] &
G \arrow[d,equal] \\
& C_2 \arrow[r,"\partial_2"'] & G .
\end{tikzcd}
\]
Equivalently, $(\beta,\mathrm{id}_G)$ is a basis-preserving isomorphism of crossed modules.
\end{pro}

\begin{proof}
Apply the universal property of $C_1$ to $C_2$ with $\alpha=\mathrm{id}_G$ and $\psi=\phi_2$. This gives a morphism $\beta\colon C_1\to C_2$. Similarly, we obtain a morphism $\gamma\colon C_2\to C_1$. By uniqueness,
$\gamma\circ\beta=\mathrm{id}_{C_1}$ and $\beta\circ\gamma=\mathrm{id}_{C_2}$.
\end{proof}

By Proposition~\ref{pro:free-CM-determined-by-f}, a free crossed module is determined by a map $f\colon S\to G$ up to a unique basis-preserving isomorphism over $\mathrm{id}_G$. If $\partial\colon C\to G$ has a basis map $\phi\colon S\to C$ with $\partial\circ\phi=f$, we call it the \emph{free crossed module on $f\colon S\to G$}. Thus, to prove existence, it is enough to construct one such crossed module for every map $f\colon S\to G$.

Whitehead gave the following construction of the free crossed module on $f\colon S\to G$ \cite{MR0030760}. Let $C_f$ be the group with generating set $S\times G$ and defining relations
\[
(t,h)^{-1}(s,g)(t,h)=(s,gh^{-1}f(t)h)
\quad(s,t\in S,\ g,h\in G).
\]
Define a right $G$-action on $C_f$ and a boundary map by
\[
(s,g)\cdot h=(s,gh),
\quad
\partial(s,g)=g^{-1}f(s)g.
\]

\begin{thm}[J.~H.~C.~Whitehead~\cite{MR0030760}]\label{thm:const-free-CM}
With the notation above, $\partial\colon C_f\to G$ is free on the map
\[
\phi\colon S\to C_f,\quad s\mapsto(s,1).
\]
Thus, $C_f\to G$ is the free crossed module on $f\colon S\to G$.
\end{thm}

A fundamental source of free crossed modules is the following theorem of Whitehead.
\begin{thm}[J.~H.~C.~Whitehead~\cite{MR0030760}]\label{thm:Whitehead}
Let $A$ be a path-connected space, and let $X$ be obtained from $A$ by attaching $2$-cells $\{e_s\}_{s\in S}$. Then the fundamental crossed module
\[
\partial\colon\pi_2(X,A)\to\pi_1(A)
\]
is a free crossed module with basis indexed by $\{e_s\}_{s\in S}$. If the attaching loop of $e_s$ represents $f(s)\in\pi_1(A)$, then this is the free crossed module on $f\colon S\to\pi_1(A)$.
\end{thm}

Proofs of Theorem~\ref{thm:Whitehead} can also be found in
Brown~\cite{MR0579818}, Brown--Higgins~\cite{MR0478150}, 
Guti\'errez--Hirschhorn~\cite{MR0816893},
Ratcliffe~\cite{MR0579810}, and Wagemann~\cite{MR4368859}.

\begin{cor}\label{cor:free-CM-topological-realization}
Every free crossed module is, up to a basis-preserving isomorphism, a fundamental crossed module as in Theorem~\ref{thm:Whitehead}.
\end{cor}

\begin{proof}
Let $\partial\colon C\to G$ be free on $\phi\colon S\to C$, and put $f=\partial\circ\phi$. Choose a path-connected CW complex $A$ with $\pi_1(A)\cong G$. For each $s\in S$, attach a $2$-cell $e_s$ along a loop representing $f(s)$, and let $X$ be the resulting CW complex. By Theorem~\ref{thm:Whitehead}, the fundamental crossed module $\pi_2(X,A)\to\pi_1(A)$ is the free crossed module on $f\colon S\to G$. The result follows from Proposition~\ref{pro:free-CM-determined-by-f}.
\end{proof}

\begin{eg}\label{eg:presentation-free-CM}
Consider a group presentation
\[
\mathscr P=\langle X\mid r_\lambda\ (\lambda\in\Lambda)\rangle,
\]
and put $F=F(X)$. Let $K(\mathscr P)$ be its presentation complex and let $K(X)$ be its $1$-skeleton. 
Let $\omega\colon\Lambda\to F$ be given by $\omega(\lambda)=r_\lambda$.
Then Theorem~\ref{thm:Whitehead} shows that the fundamental crossed module
\[
\partial\colon \pi_2(K(\mathscr P),K(X))\to F
\]
is free on $\omega$.
If we write $\partial_\omega\colon C_\omega\to F$ for the free crossed module on $\omega$, we obtain a basis-preserving isomorphism
\[
C_\omega\cong\pi_2(K(\mathscr P),K(X)).
\]
The free crossed module associated with a group presentation
plays a fundamental role in the theory of identities among
relations. Its kernel $\Ker\partial_\omega$ is the module of
identities among relations associated with the presentation~\cite{BrownHuebschmann1982}.
Since $K(X)$ is a graph, $\pi_2(K(X))=0$. The homotopy exact sequence of the pair $(K(\mathscr P),K(X))$ therefore gives
\[
\Ker\partial_\omega\cong\pi_2(K(\mathscr P)).
\]
\end{eg}

We next recall a characterization of free crossed modules due to Ratcliffe.
\begin{thm}[Ratcliffe~\cite{MR0579810}]\label{thm:Ratcliffe}
Let $\partial\colon C\to G$ be a crossed module, put $N\coloneqq\Image\,\partial$, and let $\{c_s\}_{s\in S}$ be an indexed subset of $C$. Then $\partial\colon C\to G$ is a free crossed module with basis $\{c_s\}_{s\in S}$ if and only if the following three conditions hold.
\begin{enumerate}
\item $C/[C,C]$ is a free $\mathbb Z[G/N]$-module with basis $\{c_s\bmod[C,C]\}_{s\in S}$.
\item $N$ is the normal closure of $\{\partial c_s\}_{s\in S}$ in $G$.
\item $\partial_*\colon H_2(C)\to H_2(N)$ is the zero map.
\end{enumerate}
\end{thm}

\begin{rem}
Ellis and Porter~\cite{MR0825179} proved that if
$\partial\colon C\to G$ is a projective crossed $G$-module
and $N=\Image\,\partial$, then
\[
H_2(N)\cong\Ker\,\partial\cap[C,C].
\]
Since every free crossed $G$-module is projective in
$\XMod_G$, this isomorphism applies in particular to
free crossed modules.
This result may be regarded as a crossed-module analogue
of Hopf's formula for the second homology group.
\end{rem}

We record some basic examples of free crossed modules.
\begin{eg}\label{eg:basic-free-CM1}
If $M$ is a free $G$-module with basis $S$, then the zero homomorphism $\partial\colon M\to G$ is free on the inclusion $S\hookrightarrow M$.
\end{eg}

\begin{eg}\label{eg:basic-free-CM2}
For the free group $F(S)$ on a set $S$, the identity map of $F(S)$ is free on the inclusion
$S\hookrightarrow F(S)$.
\end{eg}

\begin{eg}\label{eg:knot-group-identity-free}
Let $K\subset S^3$ be a knot, let $G_K$ be its knot group, and let $m\in G_K$ be a meridian. Then
$\mathrm{id}\colon G_K\to G_K$ is a free crossed module with basis $\{m\}$. Indeed, $G_K$ is normally generated by $m$, and $G_K^{\mathrm{ab}}\cong\mathbb Z$ is generated by the image of $m$. Moreover, $H_2(G_K)=0$~\cite{BurdeZieschangHeusener2014}*{Theorems~3.1 and~3.30}.
The result now follows by applying Theorem~\ref{thm:Ratcliffe} to $\partial=\mathrm{id}_{G_K}$.
\end{eg}

Example~\ref{eg:basic-free-CM1} is the case of a trivial boundary map. In Examples~\ref{eg:basic-free-CM2} and~\ref{eg:knot-group-identity-free}, the boundary map is the identity.
For more details on crossed modules and free crossed modules, see Brown--Higgins--Sivera~\cite{MR2841564}, Dyer~\cite{MR1279178}, Ellis~\cite{MR3971587}, Huebschmann~\cite{MR4666020}, and Wagemann~\cite{MR4368859}. For the history of crossed modules, see Huebschmann~\cites{MR4322147,MR4666020}.

\section{Augmented quandles and augmented racks}\label{sec:augmented-quandles-racks}
\subsection{Quandles and racks}
\begin{defn}
A \emph{quandle} is a set $Q$ equipped with a binary operation $(x,y)\mapsto x\ast y$ satisfying the following axioms.
\begin{enumerate}
\item $x\ast x=x\quad (x\in Q)$.
\item For each $x\in Q$, the map
$S_x\colon Q\to Q$
defined by $S_x(y)=y\ast x$ is bijective.
\item
$(x\ast y)\ast z=(x\ast z)\ast (y\ast z)\quad (x,y,z\in Q)$.
\end{enumerate}
A set equipped with a binary operation satisfying Axioms~(2) and~(3) is called a \emph{rack}.
\end{defn}

A map $\varphi\colon Q\to Q'$ between quandles is called a \emph{quandle homomorphism} if
\[
\varphi(x\ast y)=\varphi(x)\ast\varphi(y)
\quad(x,y\in Q).
\]
A rack homomorphism is defined by the same condition. A \emph{subquandle} of a quandle $Q$ is a subset $R\subset Q$ that is closed under the quandle operation and its inverse operation. 
In other words, $S_y(x),S_y^{-1}(x)\in R$ for all $x,y\in R$. A subrack is defined in the same way.

Quandles were introduced independently by Joyce and Matveev in connection with knot theory.
See \cites{MR2628474,MR0638121,MR0672410}.
For general background on quandles, see, for example, Nosaka~\cite{nosaka-book} and Bardakov--Elhamdadi--Singh~\cite{arXiv:2506.23175}.
For a brief history of racks, see Fenn--Rourke~\cite{MR1194995}.

\begin{eg}
Let $G$ be a group. The \emph{conjugation quandle} $\Conj(G)$ is the quandle whose underlying set is $G$ and whose quandle operation is
\[
g\ast h\coloneqq h^{-1}gh.
\]
More generally, a subquandle $Q\subset G$ of $\Conj(G)$ is called a \emph{conjugation subquandle}.
\end{eg}

\begin{eg}
Let $G$ be a group and let $\phi\colon G\to G$ be an automorphism. The \emph{generalized Alexander quandle} $\Alex(G,\phi)$ is the quandle whose underlying set is $G$ and whose quandle operation is
\[
g\ast h\coloneqq \phi(gh^{-1})h.
\]
If $G$ is abelian, $\Alex(G,\phi)$ is called an \emph{Alexander quandle}. Every Alexander quandle is isomorphic to a subquandle of a conjugation quandle \cite{MR4564617}.
\end{eg}

\begin{eg}[\cites{MR4175808,MR2799090,MR975077,MR1194995}]
Let $\Phi$ be a root system. Then $\Phi$ is a rack with operation
\[
\alpha\ast\beta\coloneqq r_\beta(\alpha),
\]
where
$r_\beta\colon\Phi\to\Phi$
is the reflection corresponding to the root $\beta\in\Phi$. Since
$\alpha\ast\alpha=r_\alpha(\alpha)=-\alpha$, the rack $\Phi$ is not a quandle.
\end{eg}

Let $X$ be a quandle or a rack. We introduce three groups associated with $X$. 
Let $\Aut(X)$ denote the automorphism group of $X$. 
In this paper, multiplication in $\Aut(X)$ is defined by
\[
(fg)(x)\coloneqq g(f(x))
\quad(f,g\in\Aut(X),\ x\in X).
\]
Thus $x\cdot f\coloneqq f(x)$ defines a right action on $X$. 
Let $\Inn(X)$ be the subgroup of $\Aut(X)$ generated by the maps
$S_x\colon X\to X$ $(x\in X)$.
It is called the \emph{inner automorphism group} of $X$ and is a normal
subgroup of $\Aut(X)$.
The \emph{associated group} $\As(X)$ of $X$ is defined by the presentation
\[
\As(X) \coloneqq \langle e_x \ (x \in X) \mid
e_{x \ast y} = e_y^{-1} e_x e_y \ (x,y \in X) \rangle.
\]
The canonical surjection
$\As(X) \twoheadrightarrow \Inn(X)$, $e_x \mapsto S_x$,
gives a right action of $\As(X)$ on $X$. Equivalently,
\[
x\cdot e_y=x\ast y\quad(x,y\in X).
\]

We next recall the associated quandle of a rack. Let $R$ be a rack, and let $\sim$ be the smallest rack congruence containing the relations $x\sim x\ast x$ for $x\in R$. Here, if
$x\mathbin{\bar\ast}y\coloneqq S_y^{-1}(x)$, a rack congruence means an equivalence relation compatible with both $\ast$ and $\bar\ast$. Hence
\[
[x]\ast[y]\coloneqq[x\ast y]
\]
defines a well-defined rack operation on $R/{\sim}$. Since $[x]\ast[x]=[x]$, the quotient $R/{\sim}$ is a quandle. 
We write
$Q(R)\coloneqq R/{\sim}$ and call it the \emph{associated quandle} of $R$.
See Fenn--Rourke~\cite{MR1194995}*{p.~354}.

\begin{pro}\label{pro:ass-qdl}
Let $R$ be a rack and let $Q(R)$ be its associated quandle. Then the canonical projection
$\pi\colon R\to Q(R)$
is a rack homomorphism. Moreover, for every quandle $Q$ and every rack homomorphism
$\beta\colon R\to Q$, there exists a unique quandle homomorphism
\[
\bar{\beta}\colon Q(R)\to Q,\quad\bar{\beta}\circ\pi=\beta.
\]
\end{pro}

\begin{proof}
By definition of the operation, $\pi$ is a rack homomorphism. If $\beta\colon R\to Q$ is a rack homomorphism to a quandle $Q$, then
\[
\beta(x\ast x)=\beta(x)\ast\beta(x)=\beta(x).
\]
Thus $\beta$ preserves the congruence defining $Q(R)$. Hence $\beta$ factors uniquely through $\pi$, and the induced map is a quandle homomorphism.
\end{proof}

\subsection{Augmented quandles and augmented racks}
\begin{defn}
An \emph{augmented rack} consists of a right $G$-set $X$ and a map
$\partial\colon X\to G$,
called the \emph{augmentation map}, satisfying
\[
\partial(x\cdot g)=g^{-1}(\partial x)g
\quad(x\in X,\ g\in G).
\]
An augmented rack is an \emph{augmented quandle} if, in addition,
\[
x\cdot\partial(x)=x\quad(x\in X).
\]
\end{defn}

If $\partial\colon C\to G$ is a crossed module, the first crossed module axiom makes $C$ an augmented rack over $G$. By the Peiffer identities, the induced rack operation is
\[
c_1\ast c_2
=c_1\cdot\partial(c_2)
=c_2^{-1}c_1c_2.
\]
In particular, $c\ast c=c$, so $C$ is an augmented quandle. 
Thus every crossed module is an augmented quandle, and every
augmented quandle is an augmented rack.
\begin{defn}
A \emph{morphism} of augmented racks or augmented quandles
\[
(\beta,\alpha)\colon(\partial\colon X\to G)\to(\delta\colon Y\to H)
\]
is a group homomorphism $\alpha\colon G\to H$ together with a map $\beta\colon X\to Y$ such that the following diagram commutes.
\begin{equation}\label{eq:morphism-augmented-rack-quandle}
\begin{tikzcd}
X \times G \arrow[r, ""] \arrow[d, "\beta \times \alpha"'] &
X \arrow[r, "\partial"] \arrow[d, "\beta"'] &
G \arrow[d, "\alpha"'] \\
Y \times H \arrow[r, ""] &
Y \arrow[r, "\delta"] &
H .
\end{tikzcd}
\end{equation}
\end{defn}
Equivalently,
\[
\beta(x\cdot g)=\beta(x)\cdot\alpha(g),
\quad
\delta(\beta(x))=\alpha(\partial(x)).
\]
Since the rack or quandle operation is given by $x\ast y=x\cdot\partial(y)$, these conditions imply
$\beta(x\ast y)=\beta(x)\ast\beta(y)$. Hence $\beta$ is automatically a rack homomorphism (respectively, a quandle homomorphism).

Augmented quandles were introduced by Joyce~\cites{MR2628474,MR0638121}. Augmented racks were introduced by Freyd--Yetter~\cite{MR1020583} under the name \emph{crossed $G$-sets}. 
The term ``augmented rack'' was coined by Fenn--Rourke~\cite{MR1194995}.
For a viewpoint on racks and augmented racks in terms of classifying spaces, see
Fenn--Rourke--Sanderson~\cite{MR1364012}.
For a recent homotopy-theoretic study of rack spaces of connected quandles, see Nosaka~\cite{arXiv:2610.38066}. Augmented racks have also been studied in finite group theory.
See Oda--Yoshida~\cites{MR1808345,MR2095572,MR2746028}.

\begin{eg}
Let $X$ be a rack (respectively, a quandle) and let $G$ be $\Aut(X)$ or $\Inn(X)$.
Together with the natural right action of $G$ on $X$, the canonical map
\[X\to G,\quad x\mapsto S_x\]
defines an augmented rack (respectively, an augmented quandle). 
Similarly, the canonical map
$X\to \As(X)$, $x\mapsto e_x$
is an augmented rack (respectively, an augmented quandle). 
For later use, we write
\[
\varepsilon_X\colon X\to\Inn(X),
\quad x\mapsto S_x,
\]
and call it the \emph{canonical inner augmentation} of $X$.

\end{eg}

\begin{eg}
Let $G$ be a group, and let $X\subset G$ be a subset invariant under conjugation. 
Equivalently, $X$ is a union of conjugacy classes of $G$.
Then $X$ is a subquandle of $\Conj(G)$, and the inclusion
$X\hookrightarrow G$
is an augmented quandle with respect to the conjugation action of $G$.
\end{eg}

\begin{eg}\label{eg:right-cosets-augmented-rack}
Let $G$ be a group, let $H$ be a subgroup of $G$, and let
$a\in G$. Equip the set $H\backslash G$ of right cosets
of $H$ in $G$ with the natural right $G$-action
$(Hg)\cdot k\coloneqq H(gk)$.
Consider the map
\[
\partial_a\colon H\backslash G\to G,
\quad
Hg\mapsto g^{-1}ag.
\]
This map is well-defined if and only if $a\in C_G(H)$.
Indeed, for every $h\in H$,
\[
(hg)^{-1}a(hg)=g^{-1}ag
\]
if and only if $h^{-1}ah=a$.
Moreover, whenever $\partial_a$ is well defined, we have
\[
\partial_a((Hg)\cdot k)
=k^{-1}\partial_a(Hg)k.
\]
Thus $\partial_a\colon H\backslash G\to G$
is an augmented rack if and only if $a\in C_G(H)$.
The induced rack operation is
\[
Hg\ast Hk=Hgk^{-1}ak.
\]
Furthermore, $\partial_a\colon H\backslash G\to G$
is an augmented quandle if and only if $a\in Z(H)$. 
Indeed, the augmented quandle condition
\[
(Hg)\cdot\partial_a(Hg)=Hg
\]
is equivalent to $Hag=Hg$, or equivalently, $a\in H$.
Together with $a\in C_G(H)$, this is equivalent to
$a\in Z(H)$.
\end{eg}

\begin{eg}
Let $\Phi$ be a root system and let $W$ be its Weyl group. Define a right $W$-action on $\Phi$ by
\[
\alpha\cdot w\coloneqq w^{-1}(\alpha)
\quad(\alpha\in\Phi,\ w\in W).
\]
Then the map
\[
\partial\colon\Phi\to W,
\quad
\alpha\mapsto r_\alpha
\]
is an augmented rack but not an augmented quandle. Indeed, $\alpha\cdot r_\alpha=-\alpha$, so the augmented quandle condition does not hold in general.
\end{eg}

\subsection{From augmented racks to augmented quandles}
Let $\partial\colon R\to G$ be an augmented rack. Let $\sim$ be the equivalence relation on $R$ generated by
\[
x\sim x\cdot\partial x\quad(x\in R).
\]
Put
$Q(R)\coloneqq R/{\sim}$,
and write $[x]$ for the equivalence class of $x\in R$. Define
$\partial_{\mathrm Q}\colon Q(R)\to G$ by
$\partial_{\mathrm Q}([x])\coloneqq\partial(x)$. This is well defined because
\begin{equation}\label{eq:associated-aq-augmentation-wd}
\partial(x\cdot\partial(x))
=\partial(x)^{-1}\partial(x)\partial(x)
=\partial(x).
\end{equation}
The right $G$-action on $R$ induces a well-defined right $G$-action on $Q(R)$ by
$[x]\cdot g\coloneqq[x\cdot g]$. Indeed,
\begin{equation}\label{eq:associated-aq-action-wd}
x\cdot g
\sim (x\cdot g)\cdot\partial(x\cdot g)
=(x\cdot g)\cdot(g^{-1}\partial(x)g)
=(x\cdot\partial(x))\cdot g.
\end{equation}
Moreover,
\[
[x]\cdot\partial_{\mathrm Q}([x])
=[x]\cdot\partial(x)
=[x\cdot\partial(x)]
=[x],
\]
and
\[
\partial_{\mathrm Q}([x]\cdot g)
=\partial_{\mathrm Q}([x\cdot g])
=\partial(x\cdot g)
=g^{-1}\partial(x)g
=g^{-1}\partial_{\mathrm Q}([x])g.
\]
Thus
$\partial_{\mathrm Q}\colon Q(R)\to G$
is an augmented quandle. We call it the \emph{associated augmented quandle} of $\partial\colon R\to G$.

\begin{pro}\label{pro:underlying-associated-quandle}
Let $\partial\colon R\to G$ be an augmented rack, with rack operation
$x\ast y\coloneqq x\cdot\partial(y)$.
Then the underlying quandle of the associated augmented quandle
$\partial_{\mathrm Q}\colon Q(R)\to G$
is the associated quandle of the rack $R$. In particular, the canonical projection
$\pi\colon R\to Q(R)$
is the same in both constructions.
\end{pro}

\begin{proof}
Let $\sim$ be the equivalence relation generated by
$x\sim x\cdot\partial(x)=x\ast x$. 
By Equation~\eqref{eq:associated-aq-action-wd}, the equivalence relation $\sim$
is preserved by the right $G$-action.
By Equation~\eqref{eq:associated-aq-augmentation-wd}, if $y\sim y'$, then $\partial(y)=\partial(y')$. Hence, if $x\sim x'$ and $y\sim y'$, then
\[
x\ast y
=x\cdot\partial(y)
\sim x'\cdot\partial(y)
=x'\cdot\partial(y')
=x'\ast y'.
\]
The inverse rack operation is
\[
x\mathbin{\bar\ast}y=x\cdot\partial(y)^{-1},
\]
so the same argument gives
\[
x\mathbin{\bar\ast}y\sim x'\mathbin{\bar\ast}y'.
\]
Thus $\sim$ is a rack congruence. Since $\sim$ is generated by the relations $x\sim x\ast x$, it is the smallest rack congruence containing these relations.

The underlying quandle operation of the associated augmented quandle is
\[
[x]\ast[y]
=[x]\cdot\partial_{\mathrm Q}([y])
=[x\cdot\partial(y)]
=[x\ast y].
\]
Therefore the two constructions agree.
\end{proof}

\begin{pro}\label{pro:projection}
Let $\partial\colon R\to G$ be an augmented rack, and let
$\partial_{\mathrm Q}\colon Q(R)\to G$ be its associated augmented quandle.
Let $\pi\colon R\to Q(R)$ be the canonical projection.
\begin{enumerate}
\item
$(\pi,\mathrm{id}_G)\colon(\partial\colon R\to G)\to(\partial_{\mathrm Q}\colon Q(R)\to G)$
is a morphism of augmented racks.
\item
Let $\delta\colon Q\to H$ be an augmented quandle and let
$(\beta,\alpha)\colon(\partial\colon R\to G)\to(\delta\colon Q\to H)$
be a morphism of augmented racks. Then $(\beta,\alpha)$ factors through
$\partial_{\mathrm Q}\colon Q(R)\to G$. More precisely, there exists a unique map
\[
\bar\beta\colon Q(R)\to Q
\]
such that
$
\bar\beta\circ\pi=\beta
$
and
$(\bar\beta,\alpha)\colon(\partial_{\mathrm Q}\colon Q(R)\to G)\to(\delta\colon Q\to H)$
is a morphism of augmented quandles.
\end{enumerate}
\end{pro}

\begin{proof}
(1) By definition, for every $x\in R$ and $g\in G$,
\[
\pi(x\cdot g)=[x\cdot g]=[x]\cdot g=\pi(x)\cdot g.
\]
Also, $\partial_{\mathrm Q}\circ\pi=\partial$. Hence $(\pi,\mathrm{id}_G)$ is a morphism of augmented racks.

(2) Since $Q$ is a quandle and $\beta\colon R\to Q$ is a rack homomorphism, Proposition~\ref{pro:ass-qdl} gives a unique quandle homomorphism
$\bar\beta\colon Q(R)\to Q$ satisfying $\bar\beta\circ\pi=\beta$. Since
$\alpha\circ\partial=\delta\circ\beta$, we have
\[
\alpha\circ\partial_{\mathrm Q}\circ\pi
=\alpha\circ\partial
=\delta\circ\beta
=\delta\circ\bar\beta\circ\pi.
\]
Since $\pi$ is surjective,
$\alpha\circ\partial_{\mathrm Q}=\delta\circ\bar\beta$. Moreover, for every $x\in R$ and $g\in G$,
\[
\bar\beta([x]\cdot g)
=\bar\beta([x\cdot g])
=\beta(x\cdot g)
=\beta(x)\cdot\alpha(g)
=\bar\beta([x])\cdot\alpha(g).
\]
Therefore $(\bar\beta,\alpha)$ is a morphism of augmented quandles. Uniqueness follows from the surjectivity of $\pi$.
\end{proof}

\section{Associated crossed modules and adjunctions}\label{sec:associated-crossed-modules}
\subsection{From augmented racks or augmented quandles to crossed modules}
Let $\partial\colon X\to G$ be an augmented rack or an augmented quandle. Regard $X$ as its underlying rack or quandle, and define its associated group by
\[
\As(X)\coloneqq
\langle e_x\ (x\in X)\mid
e_{x\cdot\partial(y)}=e_y^{-1}e_xe_y\ (x,y\in X)\rangle.
\]
Fenn and Rourke showed that the right $G$-action on $X$ induces a well-defined right $G$-action on $\As(X)$~\cite{MR1194995}*{pp.~355--356}, given by
\[
e_x\cdot g\coloneqq e_{x\cdot g}.
\]
The assignment $e_x\mapsto\partial(x)$ extends to a well-defined group homomorphism
\[
\partial_{\mathrm{As}}\colon\As(X)\to G.
\]
Together with the right $G$-action above, this homomorphism defines a crossed module. We call
$\partial_{\mathrm{As}}\colon\As(X)\to G$
the \emph{associated crossed module} of $\partial\colon X\to G$.

\begin{rem}
The associated crossed module of an augmented rack also appears in
Crans--Wagemann~\cite{CransWagemann2014}*{Proposition~2.23}
in the context of crossed modules of racks.
\end{rem}

\begin{pro}\label{pro:As-associated-quandle}
Let $\partial\colon R\to G$ be an augmented rack.
Let $\partial_{\mathrm Q}\colon Q(R)\to G$ be its associated augmented quandle.
The canonical projection $R\to Q(R)$ induces an isomorphism
\[
\As(R)\xrightarrow{\cong}\As(Q(R)),
\quad e_x\mapsto e_{[x]}.
\]
Moreover, this is an isomorphism of crossed modules over $\mathrm{id}_G$.
\end{pro}

\begin{proof}
Since the canonical projection
$\pi\colon R\to Q(R)$, $x\mapsto[x]$,
is a rack homomorphism, it induces a group homomorphism
\[
\pi_*\colon\As(R)\to\As(Q(R)),
\quad e_x\mapsto e_{[x]}.
\]
We construct an inverse. The defining relation of $\As(R)$ gives
\[
e_{x\ast x}=e_x^{-1}e_xe_x=e_x
\quad(x\in R).
\]
The map
$R\to\Conj(\As(R))$, $x\mapsto e_x$,
is a rack homomorphism because
\[
e_{x\ast y}=e_y^{-1}e_xe_y.
\]
It also sends $x$ and $x\ast x$ to the same element. By the universal property in Proposition~\ref{pro:ass-qdl}, this map factors uniquely through
$\pi\colon R\to Q(R)$. Hence we obtain a quandle homomorphism
\[
Q(R)\to\Conj(\As(R)),
\quad [x]\mapsto e_x.
\]
By the defining presentation of the associated group, this induces a group homomorphism
\[
\rho\colon\As(Q(R))\to\As(R),
\quad e_{[x]}\mapsto e_x.
\]
The compositions $\rho\circ\pi_*$ and $\pi_*\circ\rho$ are the identity on the generators $e_x$ and $e_{[x]}$, respectively. Therefore $\pi_*$ and $\rho$ are inverse to each other, and
\[
\As(R)\cong\As(Q(R)).
\]

It remains to check that this is an isomorphism of crossed modules. For the boundary maps,
\[
\partial_{\mathrm{As}(Q(R))}(\pi_*(e_x))
=\partial_{\mathrm{As}(Q(R))}(e_{[x]})
=\partial_{\mathrm Q}([x])
=\partial(x)
=\partial_{\mathrm{As}(R)}(e_x).
\]
Thus $\pi_*$ is compatible with the boundary maps. Also, for every $g\in G$,
\[
\pi_*(e_x\cdot g)
=\pi_*(e_{x\cdot g})
=e_{[x\cdot g]}
=e_{[x]\cdot g}
=e_{[x]}\cdot g
=\pi_*(e_x)\cdot g.
\]
Hence $\pi_*$ is compatible with the $G$-actions. Therefore $\pi_*$ is an isomorphism of crossed modules over $\mathrm{id}_G$.
\end{proof}

\subsection{Functoriality and adjunctions}
We first fix notation.
Let $\XMod$, $\AugQdl$, and $\AugRack$ denote the categories of crossed modules, augmented quandles, and augmented racks, respectively.
A crossed module is an augmented quandle, and an augmented quandle is
an augmented rack. We write the corresponding forgetful functors as
\[
U_{\mathrm Q}\colon\XMod\to\AugQdl,
\quad
U_{\mathrm R}\colon\AugQdl\to\AugRack.
\]
We now formulate the associated constructions as left adjoints to these
forgetful functors.

\begin{thm}\label{thm:adjunctions}
The associated constructions define functors
\[
Q\colon\AugRack\to\AugQdl,
\quad
\As\colon\AugQdl\to\XMod,
\]
and there are adjunctions
\[
Q\dashv U_{\mathrm R},
\quad
\As\dashv U_{\mathrm Q}.
\]
Consequently,
$\As\circ Q\dashv U_{\mathrm R}\circ U_{\mathrm Q}$.
Moreover, by Proposition~\ref{pro:As-associated-quandle}, the functor
\[
\As\circ Q\colon\AugRack\to\XMod
\]
is naturally isomorphic to the associated crossed module functor defined directly on augmented racks.
\end{thm}

\begin{proof}
We first describe the associated constructions on morphisms. Let
$(\beta,\alpha)\colon(R\to G)\to(R'\to H)$ be a morphism of augmented racks, and define
$Q(\beta)([x])\coloneqq[\beta(x)]$.
Then $(Q(\beta),\alpha)$ is a morphism of the associated augmented quandles. Similarly, for a morphism of augmented quandles
$(\beta,\alpha)\colon(X\to G)\to(Y\to H)$, define
$\As(\beta)(e_x)\coloneqq e_{\beta(x)}$ for $x\in X$. Then $(\As(\beta),\alpha)$ is a morphism of the associated crossed modules. These constructions preserve identities and compositions. Hence $Q$ and $\As$ define functors.

By the universal property in Proposition~\ref{pro:projection}, for every augmented rack $R$ and every augmented quandle $Q$, there is a natural bijection
\[
\operatorname{Hom}_{\AugQdl}(Q(R),Q)
\cong
\operatorname{Hom}_{\AugRack}(R,U_{\mathrm R}(Q)).
\]
Thus $Q\dashv U_{\mathrm R}$.

We prove the second adjunction $\As\dashv U_{\mathrm Q}$. Let
$\partial\colon X\to G$ be an augmented quandle, let
$\delta\colon C\to H$ be a crossed module, and let
\[
(\beta,\alpha)\colon
(\partial\colon X\to G)
\to
U_{\mathrm Q}(\delta\colon C\to H)
\]
be a morphism of augmented quandles. Consider the assignment
$e_x\mapsto\beta(x)$ for $x\in X$. By the Peiffer identities, the underlying quandle operation on $C$ is
\[
c_1\ast c_2
=c_1\cdot\delta(c_2)
=c_2^{-1}c_1c_2.
\]
Therefore
\[
\begin{aligned}
\beta(x\ast y)
&=\beta(x\cdot\partial(y))=\beta(x)\cdot\alpha(\partial(y))\\
&=\beta(x)\cdot\delta(\beta(y))=\beta(y)^{-1}\beta(x)\beta(y).
\end{aligned}
\]
Thus the defining relations of $\As(X)$ are preserved, and the assignment
$e_x\mapsto\beta(x)$ extends uniquely to a group homomorphism
$\bar\beta\colon\As(X)\to C$.
For each $x\in X$,
\[
\delta(\bar\beta(e_x))
=\delta(\beta(x))
=\alpha(\partial(x))
=\alpha(\partial_{\mathrm{As}}(e_x)).
\]
For every $g\in G$,
\[
\begin{aligned}
\bar\beta(e_x\cdot g)
&=\bar\beta(e_{x\cdot g})=\beta(x\cdot g)\\
&=\beta(x)\cdot\alpha(g)=\bar\beta(e_x)\cdot\alpha(g).
\end{aligned}
\]
Since the elements $e_x$ $(x\in X)$ generate $\As(X)$, the pair $(\bar\beta,\alpha)$ is a morphism of crossed modules
\[
(\bar\beta,\alpha)\colon
(\partial_{\mathrm{As}}\colon\As(X)\to G)
\to
(\delta\colon C\to H).
\]

Conversely, let
\[
(\bar\beta,\alpha)\colon
(\partial_{\mathrm{As}}\colon\As(X)\to G)
\to
(\delta\colon C\to H)
\]
be a morphism of crossed modules. Define
\[
\beta\colon X\to C,
\quad
\beta(x)\coloneqq\bar\beta(e_x).
\]
Then
\[
\delta(\beta(x))
=\delta(\bar\beta(e_x))
=\alpha(\partial_{\mathrm{As}}(e_x))
=\alpha(\partial(x)),
\]
and for every $g\in G$,
\[
\begin{aligned}
\beta(x\cdot g)
&=\bar\beta(e_{x\cdot g})=\bar\beta(e_x\cdot g)\\
&=\bar\beta(e_x)\cdot\alpha(g)=\beta(x)\cdot\alpha(g).
\end{aligned}
\]
Hence $(\beta,\alpha)$ is a morphism of augmented quandles
\[
(\partial\colon X\to G)
\to
U_{\mathrm Q}(\delta\colon C\to H).
\]

These two constructions are inverse to each other by the uniqueness of the extension from the generators $e_x$. Hence there is a natural bijection
\[
\operatorname{Hom}_{\XMod}(\As(X),C)
\cong
\operatorname{Hom}_{\AugQdl}(X,U_{\mathrm Q}(C)).
\]
Therefore $\As\dashv U_{\mathrm Q}$.
Finally, Proposition~\ref{pro:As-associated-quandle} gives a natural isomorphism between $\As\circ Q$ and the associated crossed module functor defined directly on augmented racks. Naturality follows immediately from the formula $e_x\mapsto e_{[x]}$.
\end{proof}

\section{Free augmented racks, quandles, and crossed modules}\label{sec:free-augmented}
Following the definition of a free crossed module, we introduce free augmented racks and free augmented quandles by analogous universal properties. Let $f\colon S\to G$ be a map. We construct the corresponding free augmented rack and free augmented quandle and compare their associated crossed modules with the free crossed module on $f$.

\begin{defn}
Let $\partial\colon X\to G$ be an augmented rack (respectively, an augmented quandle), let $S$ be a set, and let $\phi\colon S\to X$ be a map. We say that $\partial\colon X\to G$ is \emph{free on $\phi$} if the following universal property holds. For every augmented rack (respectively, augmented quandle) $\delta\colon Y\to H$, every group homomorphism $\alpha\colon G\to H$, and every map $\psi\colon S\to Y$ satisfying $\alpha\circ\partial\circ\phi=\delta\circ\psi$, there exists a unique map $\beta\colon X\to Y$ such that $\beta\circ\phi=\psi$ and $(\beta,\alpha)$ is a morphism of augmented racks (respectively, augmented quandles).
\[
\begin{tikzcd}
S \arrow[r,"\phi"] \arrow[rd,"\psi"'] &
X \arrow[r,"\partial"] \arrow[d,"\exists!\,\beta",dotted] &
G \arrow[d,"\alpha"] \\
& Y \arrow[r,"\delta"'] & H .
\end{tikzcd}
\]
\end{defn}

In this situation, $\phi$ is injective. Indeed, give $S$ the trivial rack (respectively, trivial quandle) structure and consider the augmented rack (respectively, augmented quandle) $S\to1$. Apply the universal property to the unique homomorphism $G\to1$ and to $\mathrm{id}_S\colon S\to S$. This gives a left inverse $X\to S$ of $\phi$, so $\phi$ is injective. We identify $S$ with its image and call $S$ a \emph{basis} of the free augmented rack (respectively, free augmented quandle).

\begin{rem}
We allow the empty set to be a rack and a quandle.
We also allow the basis of a free augmented rack, free augmented quandle,
or free crossed module to be empty.
Let $f\colon \varnothing\to G$ be the unique map.
Then the free augmented rack and the free augmented quandle on $f$
are both
$
\varnothing\to G$,
whereas the free crossed module on $f$ is
$
1\to G$.
\end{rem}

As in the case of free crossed modules, the free object is determined by the composite
$f=\partial\circ\phi\colon S\to G$.

\begin{pro}\label{pro:free-augmented-determined-by-f}
Let $G$ be a group, let $S$ be a set, and let $f\colon S\to G$ be a map. For $i=1,2$, suppose that $\partial_i\colon X_i\to G$ is free on $\phi_i\colon S\to X_i$ as an augmented rack (respectively, augmented quandle), and that $\partial_i\circ\phi_i=f$. Then there exists a unique isomorphism
\[
(\beta,\mathrm{id}_G)\colon(\partial_1\colon X_1\to G)
\xrightarrow{\cong}
(\partial_2\colon X_2\to G)
\]
such that $\beta\circ\phi_1=\phi_2$.
\end{pro}

\begin{proof}
Apply the universal property of $X_1$ to $X_2$ with $\alpha=\mathrm{id}_G$ and $\psi=\phi_2$. This gives a morphism $\beta\colon X_1\to X_2$. Similarly, we obtain a morphism $\gamma\colon X_2\to X_1$. By uniqueness,
$\gamma\circ\beta=\mathrm{id}_{X_1}$ and $\beta\circ\gamma=\mathrm{id}_{X_2}$.
\end{proof}

If $\partial\colon X\to G$ has a basis map $\phi\colon S\to X$ with $\partial\circ\phi=f$, we call it the \emph{free augmented rack on $f\colon S\to G$} (respectively, the \emph{free augmented quandle on $f\colon S\to G$}).

\begin{thm}\label{thm:AT1}
Let $G$ be a group, let $S$ be a set, and let $f\colon S\to G$ be a map. Define a right $G$-action on $S\times G$ by
\[
(s,g)\cdot h\coloneqq(s,gh).
\]
Then
\[
\partial_f\colon S\times G\to G,
\quad
\partial_f(s,g)=g^{-1}f(s)g,
\]
is the free augmented rack on $f$ with basis $s\mapsto(s,1)$. We denote this free augmented rack by $R_f\to G$.
\end{thm}

\begin{proof}
First,
\[
\partial_f((s,g)\cdot h)
=h^{-1}\partial_f(s,g)h,
\]
so $\partial_f$ is an augmented rack. Let $\delta\colon R\to H$ be an augmented rack, let $\alpha\colon G\to H$ be a group homomorphism, and let $\psi\colon S\to R$ satisfy
$\alpha\circ f=\delta\circ\psi$. Put $x_s=\psi(s)$ and define
\[
\beta(s,g)=x_s\cdot\alpha(g).
\]
Then
\[
\delta(\beta(s,g))
=\alpha(g)^{-1}\delta(x_s)\alpha(g)
=\alpha(g^{-1}f(s)g)
=\alpha(\partial_f(s,g)),
\]
and
$\beta((s,g)\cdot h)=\beta(s,g)\cdot\alpha(h)$. Hence $(\beta,\alpha)$ is a morphism of augmented racks and $\beta(s,1)=\psi(s)$.

Uniqueness follows from equivariance. Indeed, every element of $S\times G$ has the form $(s,g)=(s,1)\cdot g$. Thus, if $\beta'$ satisfies the same conditions, then
\[
\beta'(s,g)
=\beta'(s,1)\cdot\alpha(g)
=\psi(s)\cdot\alpha(g)
=\beta(s,g).
\qedhere\]
\end{proof}

\begin{thm}\label{thm:AT2}
Let $R_f=S\times G\to G$ be the free augmented rack as in Theorem~\ref{thm:AT1}. 
Its associated augmented quandle
\[
\partial_{\mathrm Q}\colon Q(S\times G)\to G,
\quad
[s,g]\mapsto g^{-1}f(s)g,
\]
is the free augmented quandle on $f$ with basis $s\mapsto[s,1]$. We denote this free augmented quandle by $Q_f\to G$.
\end{thm}

\begin{proof}
Let $\delta\colon Q\to H$ be an augmented quandle, let $\alpha\colon G\to H$ be a group homomorphism, and let $\psi\colon S\to Q$ satisfy
$\alpha\circ f=\delta\circ\psi$. Regard $Q$ as its underlying augmented rack. By Theorem~\ref{thm:AT1}, there exists a unique morphism of augmented racks
\[
(\beta,\alpha)\colon
(R_f\to G)\to(Q\to H)
\]
extending $\psi$. Thus
$\beta(s,1)=\psi(s)$ for all $s\in S$.

Let $\pi\colon R_f\to Q_f=Q(R_f)$ be the canonical projection. By Proposition~\ref{pro:projection}, $(\beta,\alpha)$ factors uniquely through $\pi$. Therefore there exists a unique morphism of augmented quandles
\[
(\bar\beta,\alpha)\colon
(Q_f\to G)\to(Q\to H)
\]
such that
$\bar\beta\circ\pi=\beta$. In particular,
\[
\bar\beta([s,1])
=\beta(s,1)
=\psi(s)
\quad(s\in S).
\]

We verify uniqueness. Suppose that $\gamma\colon Q_f\to Q$ is such that $(\gamma,\alpha)$ is a morphism of augmented quandles and
$\gamma([s,1])=\psi(s)$ for every $s\in S$. Then
\[
\gamma\circ\pi\colon R_f\to Q
\]
is a morphism of augmented racks extending $\psi$. By the uniqueness in Theorem~\ref{thm:AT1},
\[
\gamma\circ\pi=\beta
=\bar\beta\circ\pi.
\]
Since $\pi$ is surjective, $\gamma=\bar\beta$. Thus $Q_f$ is the free augmented quandle on $f$.
\end{proof}

\begin{cor}\label{cor:free-augmented-existence-uniqueness}
For every map $f\colon S\to G$, the free augmented rack $R_f\to G$ and the free augmented quandle $Q_f\to G$ exist. Each is unique up to a unique basis-preserving isomorphism over $\mathrm{id}_G$.
\end{cor}

\begin{thm}\label{thm:AT3}
Let $f\colon S\to G$ be a map. The associated crossed modules of $R_f\to G$ and $Q_f\to G$ are free crossed modules on $f$. Let $C_f\to G$ denote the free crossed module determined by $f$, and let $c_s\in C_f$ be the basis element corresponding to $s\in S$. Then there exist unique isomorphisms of crossed modules over $\mathrm{id}_G$
\[
\begin{aligned}
\Theta_R&\colon\As(R_f)\xrightarrow{\cong}C_f,
& e_{(s,1)}&\mapsto c_s,\\
\Theta_Q&\colon\As(Q_f)\xrightarrow{\cong}C_f,
& e_{[s,1]}&\mapsto c_s.
\end{aligned}
\]
\end{thm}

\begin{proof}
Let $\delta\colon D\to H$ be a crossed module, let $\alpha\colon G\to H$ be a group homomorphism, and let $\psi\colon S\to D$ satisfy
$\alpha\circ f=\delta\circ\psi$. Regard $D\to H$ as its underlying augmented quandle. By the freeness of $Q_f$, there exists a unique morphism of augmented quandles
\[
(\beta,\alpha)\colon
(Q_f\to G)\to(D\to H)
\]
extending $\psi$. In particular,
\[
\beta([s,1])=\psi(s)
\quad(s\in S).
\]
By the adjunction $\As\dashv U_{\mathrm Q}$, this corresponds to a unique morphism of crossed modules
\[
(\bar\beta,\alpha)\colon
(\As(Q_f)\to G)\to(D\to H)
\]
such that
\[
\bar\beta(e_{[s,1]})=\beta([s,1])=\psi(s)
\quad(s\in S).
\]
Thus $\As(Q_f)\to G$ is the free crossed module on $f$ with basis $\{e_{[s,1]}\}_{s\in S}$.

By Proposition~\ref{pro:As-associated-quandle}, $\As(R_f)\to G$ is also the free crossed module on $f$, with basis $\{e_{(s,1)}\}_{s\in S}$. Proposition~\ref{pro:free-CM-determined-by-f} then gives the stated unique isomorphisms sending these basis elements to $c_s$ $(s\in S)$.
\end{proof}

\begin{cor}\label{cor:free-core}
For every map $f\colon S\to G$, there exist basis-preserving isomorphisms over $\mathrm{id}_G$
\[
\As(R_f)\cong\As(Q_f)\cong C_f.
\]
Consequently, every free crossed module is the associated crossed module of both a free augmented rack and a free augmented quandle. Moreover, every free augmented quandle can be realized as the associated augmented quandle of a free augmented rack.
\end{cor}

\begin{proof}
The first isomorphisms follow from Theorem~\ref{thm:AT3}. A free crossed module with basis map $\phi$ and boundary map $\partial$ is isomorphic to $C_f$, where $f=\partial\circ\phi$. Similarly, every free augmented quandle is isomorphic to $Q_f=Q(R_f)$ for the map $f$ determined by its basis and augmentation map. In each case, the isomorphism preserves the basis.
\end{proof}

\begin{rem}
The relations $[s,g]=[s,f(s)g]$ defining $Q_f$ give no new relations after taking the associated group. Indeed, in the associated group of $R_f$, put $x=(s,g)$. Then
\[
e_{x\cdot\partial(x)}=e_x^{-1}e_xe_x=e_x.
\]
Since $x\cdot\partial(x)=(s,f(s)g)$, we obtain
$e_{(s,f(s)g)}=e_{(s,g)}$. 
\end{rem}

\subsection{Free racks and free quandles}
Let $G=F(S)$ and let $f\colon S\hookrightarrow F(S)$ be the canonical inclusion. Then the constructions above recover the standard models of the classical free rack and the classical free quandle.
See Fenn--Rourke~\cite{MR1194995} and Joyce~\cite{MR0638121} for these models.
The \emph{free rack} $\mathrm{FR}(S)$ on $S$ is the set $S\times F(S)$ with operation
\[
(s,g)\ast(t,h)=(s,gh^{-1}th).
\]
The \emph{free quandle} $\mathrm{FQ}(S)$ on $S$ is the quotient of $\mathrm{FR}(S)$ by the equivalence relation generated by
\[
(s,g)\sim(s,sg),
\]
with quandle operation
\[
[s,g]\ast[t,h]=[s,gh^{-1}th].
\]

\begin{pro}\label{pro:free-rack-quandle-group}
Let $S$ be a set and let $F(S)$ be the free group on $S$. Consider the free augmented rack $S\times F(S)\to F(S)$ corresponding to the inclusion $S\hookrightarrow F(S)$, its associated augmented quandle $Q(S\times F(S))\to F(S)$, and the associated crossed module $\As(S\times F(S))\to F(S)$. Then:
\begin{enumerate}
\item $S\times F(S)=\mathrm{FR}(S)$ as racks.
\item $Q(S\times F(S))\cong\mathrm{FQ}(S)$ as quandles.
\item $\As(S\times F(S))\cong F(S)$ as crossed modules over $F(S)$.
\end{enumerate}
\end{pro}

\begin{proof}
Apply Theorem~\ref{thm:AT1} with $G=F(S)$ and the inclusion $S\hookrightarrow F(S)$. The underlying rack operation is
\[
(s,g)\ast(t,h)=(s,gh^{-1}th),
\]
which proves (1). By Proposition~\ref{pro:underlying-associated-quandle}, the underlying quandle of the associated augmented quandle is the associated quandle of $\mathrm{FR}(S)$, namely $\mathrm{FQ}(S)$. This proves (2). Finally, Theorem~\ref{thm:AT3} shows that $\As(S\times F(S))\to F(S)$ is the free crossed module on the inclusion $S\hookrightarrow F(S)$. Hence (3) follows from Example~\ref{eg:basic-free-CM2} and Proposition~\ref{pro:free-CM-determined-by-f}.
\end{proof}

\begin{rem}\label{rem:free-quandle-two-models}
For the inclusion $\iota\colon S\hookrightarrow F(S)$, the associated free augmented quandle is
\[
Q_\iota=(S\times F(S))/\sim,
\]
where
\[
(s,g)\sim(t,h)
\quad\Longleftrightarrow\quad
s=t\ \text{and}\ h=s^n g\ \text{for some }n\in\mathbb Z.
\]
Hence
\[
Q_\iota\cong\bigsqcup_{s\in S}\langle s\rangle\backslash F(S).
\]
The augmentation
\[
Q_\iota\to F(S),
\quad
\langle s\rangle g\mapsto g^{-1}sg,
\]
identifies $Q_\iota$ with a subquandle of $\Conj(F(S))$.
Its image consists of all conjugates of the free generators of $F(S)$. 
Indeed, distinct free generators are not conjugate, and
$C_{F(S)}(s)=\langle s\rangle$ for each $s\in S$, so this map is injective. 
This is the standard isomorphism between Joyce's coset model and the conjugation model of the free quandle \cite{MR0638121}.
\end{rem}

\section{Stabilizers and a freeness criterion}\label{sec:stabilizers}
The description of quandles by cosets goes back to Joyce~\cite{MR2628474}. 
Here we describe the right $G$-action on an augmented rack or augmented quandle
in terms of right cosets of stabilizers.
For an augmented quandle, the condition
\[
x\cdot\partial(x)=x
\]
implies that the stabilizer of $x$ contains $\langle\partial(x)\rangle$. We show below that freeness is characterized by the case in which this stabilizer is minimal, namely
\[
\Stab_G(x)=\langle\partial(x)\rangle.
\]

\subsection{Stabilizer quotients}

\begin{pro}\label{pro:stabilizer-quotient}
Let $\partial\colon X\to G$ be an augmented rack, and let $S\subset X$
be a complete set of representatives for the $G$-orbits. Put
\[
f=\partial|_S\colon S\to G,
\quad
H_s=\Stab_G(s)\quad(s\in S),
\]
and
\[
Y=\bigsqcup_{s\in S}H_s\backslash G.
\]
Then
$H_s\le C_G(f(s))$ $(s\in S)$.
Define a right $G$-action on $Y$ by
$(H_sg)\cdot k=H_sgk$,
and define
\[
\partial_Y\colon Y\to G,
\quad
\partial_Y(H_sg)=g^{-1}f(s)g.
\]
Then $\partial_Y\colon Y\to G$ is an augmented rack, and
\[
\Theta\colon Y\xrightarrow{\cong}X,
\quad
H_sg\mapsto s\cdot g
\]
is an isomorphism of augmented racks.
The rack operation on $Y$ is
\[
(H_sg)\ast(H_th)=H_sgh^{-1}f(t)h.
\]
\end{pro}

\begin{proof}
If $a\in H_s$, then $s\cdot a=s$, and hence
\[
f(s)=\partial(s)
=\partial(s\cdot a)
=a^{-1}f(s)a.
\]
Thus $H_s\le C_G(f(s))$, and consequently $\partial_Y$ is well-defined.
Moreover,
\[
\partial_Y((H_sg)\cdot k)
=k^{-1}\partial_Y(H_sg)k,
\]
so $\partial_Y\colon Y\to G$ is an augmented rack.

Since $S$ is a complete set of representatives for the $G$-orbits,
$\Theta$ is surjective. If
$s\cdot g=t\cdot h$,
then $s=t$, and for a fixed $s$,
\[
s\cdot g=s\cdot h
\quad\Longleftrightarrow\quad
gh^{-1}\in H_s
\quad\Longleftrightarrow\quad
H_sg=H_sh.
\]
Thus $\Theta$ is bijective. It is compatible with the right $G$-actions,
and
\[
\partial(\Theta(H_sg))
=\partial(s\cdot g)
=g^{-1}f(s)g
=\partial_Y(H_sg).
\]
Hence $\Theta$ is an isomorphism of augmented racks. Finally,
\[
(H_sg)\ast(H_th)
=(H_sg)\cdot\partial_Y(H_th)
=H_sgh^{-1}f(t)h.
\qedhere\]
\end{proof}

\begin{cor}\label{cor:stabilizer-canonical-morphisms}
With the notation of Proposition~\ref{pro:stabilizer-quotient},
the unique morphism of augmented racks
\[
(\eta,\mathrm{id}_G)\colon
(R_f\to G)\to(\partial\colon X\to G)
\]
extending the inclusion $S\hookrightarrow X$ is, under the identification
$X\cong Y$, given by
\[
S\times G\to Y,
\quad
(s,g)\mapsto H_sg.
\]
If $\partial\colon X\to G$ is an augmented quandle, then
\[
\langle f(s)\rangle\le H_s
\quad(s\in S),
\]
and $\Theta\colon Y\to X$ is an isomorphism of augmented quandles.
Moreover, under the identification
\[
Q_f\cong\bigsqcup_{s\in S}\langle f(s)\rangle\backslash G,
\]
the unique morphism of augmented quandles
\[
(\bar\eta,\mathrm{id}_G)\colon
(Q_f\to G)\to(\partial\colon X\to G)
\]
extending $S\hookrightarrow X$ is given by
$\langle f(s)\rangle g\mapsto H_sg$.
\end{cor}

\begin{proof}
Since $\eta$ is $G$-equivariant and $\eta(s,1)=s$, we have
$\eta(s,g)=s\cdot g$.
Under the identification $X\cong Y$, this is the map
\[
(s,g)\mapsto H_sg.
\]

Suppose that $X\to G$ is an augmented quandle. Since $s\cdot f(s)=s$, we have
$f(s)\in H_s$, and hence $\langle f(s)\rangle\le H_s$.
It follows that $Y\to G$ is an augmented quandle and that $\Theta$ is an
isomorphism of augmented quandles. The map from $R_f$ therefore factors
through $Q_f$, and the induced map is
\[
\langle f(s)\rangle g\mapsto H_sg.
\qedhere\]
\end{proof}

\begin{cor}\label{cor:injective-augmentation-centralizer}
Let $\partial\colon Q\to G$ be an injective augmented quandle, and let
$S\subset Q$ be a complete set of representatives for the $G$-orbits.
Then
\[
\Stab_G(s)=C_G(\partial(s))
\quad(s\in S).
\]
Consequently,
\[
Q\cong
\bigsqcup_{s\in S}C_G(\partial(s))\backslash G,
\quad
C_G(\partial(s))g\mapsto s\cdot g.
\]
Moreover, $\partial$ identifies $Q$ with a union of conjugacy classes in
$G$, viewed as a subquandle of $\Conj(G)$.
\end{cor}

\begin{proof}
By Proposition~\ref{pro:stabilizer-quotient},
$\Stab_G(s)\le C_G(\partial(s))$.
Conversely, if $a\in C_G(\partial(s))$, then
\[
\partial(s\cdot a)
=a^{-1}\partial(s)a
=\partial(s).
\]
Since $\partial$ is injective, $s\cdot a=s$, and hence
$a\in\Stab_G(s)$. Thus
\[
\Stab_G(s)=C_G(\partial(s)).
\]
The remaining statements follow from
Proposition~\ref{pro:stabilizer-quotient}.
\end{proof}

\begin{rem}\label{rem:GA-quandle}
Let $A\subset G$ be a nonempty subset, and let
$f\colon A\hookrightarrow G$ be the inclusion. The $(G,A)$-quandle
$Q(G,A)$ of Bardakov--Nasybullov~\cite{BardakovNasybullov2019}*{Section~4}
is given by
\[
Q(G,A)\cong\bigsqcup_{a\in A}C_G(a)\backslash G.
\]
The canonical morphism from the free augmented quandle
\[
Q_f\cong\bigsqcup_{a\in A}\langle a\rangle\backslash G
\]
to $Q(G,A)$ is given by
\[
\langle a\rangle g\mapsto C_G(a)g.
\]
If the elements of $A$ are pairwise non-conjugate, then
\[
C_G(a)g\mapsto g^{-1}ag
\]
gives an isomorphism from $Q(G,A)$ to the union of the conjugacy classes
of the elements of $A$ in $\Conj(G)$
\cite{BardakovNasybullov2019}*{Proposition~4.3}.
\end{rem}

\subsection{A freeness criterion}
By Corollary~\ref{cor:free-core}, the associated crossed module of a free augmented quandle is free. The converse does not hold in general. The stabilizer description above gives the following characterization of freeness for augmented quandles.

\begin{pro}\label{pro:free-aq-coset-criterion}
Let $\partial\colon Q\to G$ be an augmented quandle, let $i\colon S\to Q$ be a map, and put
$f=\partial\circ i\colon S\to G$. Assume that $Q=i(S)\cdot G$. Then there is a surjective morphism of augmented quandles
\[
\Phi\colon\bigsqcup_{s\in S}\langle f(s)\rangle\backslash G\to Q,
\quad
\langle f(s)\rangle g\mapsto i(s)\cdot g.
\]
Moreover, $\partial\colon Q\to G$ is a free augmented quandle with basis $i\colon S\to Q$ if and only if $\Phi$ is an isomorphism.
In particular, suppose that $i\colon S\to Q$ is injective and that $i(S)$ is a complete set of representatives for the $G$-orbits. Then $\Phi$ is an isomorphism if and only if
\[
\Stab_G(i(s))=\langle f(s)\rangle
\quad(s\in S).
\]
\end{pro}

\begin{proof}
By the construction in Theorem~\ref{thm:AT2}, the free augmented quandle on $f$ is
\[
Q_f\cong\bigsqcup_{s\in S}\langle f(s)\rangle\backslash G.
\]
By the universal property, there is a unique morphism
$\Phi\colon Q_f\to Q$ sending the basis element $\langle f(s)\rangle$ to $i(s)$. The assumption $Q=i(S)\cdot G$ implies that $\Phi$ is surjective.
If $\Phi$ is an isomorphism, then $Q$ is free on $i$. Conversely, if
$Q$ is free on $i$, Proposition~\ref{pro:free-augmented-determined-by-f}
shows that $\Phi$ is the unique basis-preserving isomorphism from $Q_f$
to $Q$. Hence $Q$ is free on $i$ if and only if $\Phi$ is an isomorphism.

Now assume that $i$ is injective and that $i(S)$ is a complete set of orbit representatives. Put
$H_s=\Stab_G(i(s))$. By Proposition~\ref{pro:stabilizer-quotient}, on each orbit the map $\Phi$ is identified with
\[
\langle f(s)\rangle\backslash G\to H_s\backslash G,
\quad
\langle f(s)\rangle g\mapsto H_sg.
\]
Therefore $\Phi$ is an isomorphism if and only if
$H_s=\langle f(s)\rangle$ for every $s\in S$.
\end{proof}

\section{Examples and applications of free augmented quandles}\label{sec:examples-applications}
We apply the stabilizer criterion to several classes of augmented
quandles. After some basic examples, we consider conjugation quandles,
generalized Alexander quandles, and knot quandles.

\subsection{The constant map to the identity}
Let $S$ be a set and let $G$ be a group. Define $f\colon S\to G$ by $f(s)=1$ $(s\in S)$. Then
\[
R_f=S\times G,
\qquad
\partial_f(s,g)=1.
\]
Thus $R_f$ already satisfies the augmented quandle condition, and the projection to the associated augmented quandle is an isomorphism. Hence we may identify
\[
Q_f=R_f=S\times G.
\]
The underlying rack or quandle operation is trivial:
\[
(s,g)\ast(t,h)=(s,g).
\]
The associated crossed module is the zero crossed module
\[
0\colon M_S\to G,
\]
where $M_S$ is the free $G$-module with basis $S$. This follows from Example~\ref{eg:basic-free-CM1} and Theorem~\ref{thm:AT3}.

\subsection{A singleton basis}
Let $S=\{s\}$ and let $f(s)=a\in G$. Then the free augmented rack
$R_f=\{s\}\times G$ may be identified with $G$, and its augmentation is
\[
\partial_f(g)=g^{-1}ag.
\]
The free augmented quandle is
\[
Q_f\cong\langle a\rangle\backslash G,
\]
with augmentation and quandle operation
\[
\partial_f(\langle a\rangle g)=g^{-1}ag,
\qquad
\langle a\rangle g\ast\langle a\rangle h
=\langle a\rangle gh^{-1}ah.
\]

By Theorem~\ref{thm:AT3}, the associated crossed modules of $R_f$ and $Q_f$
are isomorphic to the free crossed module
\[
C_a\to G
\]
with a one-element basis $\{c\}$ such that $\partial(c)=a$.
Equivalently, $C_a\to G$ may be realized as the fundamental crossed module
of the pair
\[
\bigl(K(G,1)\cup_a e^2,K(G,1)\bigr).
\]
In particular, when $a=1$, this is the zero crossed module on the free $G$-module with one generator from the preceding subsection.

Suppose further that $C_G(a)=\langle a\rangle$. Then the map
$
\langle a\rangle g\mapsto g^{-1}ag
$
is a bijection from $\langle a\rangle\backslash G$ onto the conjugacy class of $a$. Thus $Q_f$ may be identified with the conjugacy class of $a$.

\subsection{Conjugation quandles}
Let $S\subset G$ and put
\[
Q=\bigcup_{g\in G}g^{-1}Sg\subset\Conj(G). 
\]
Then the inclusion
$\partial\colon Q\hookrightarrow G$
is an augmented quandle. Corollary~\ref{cor:injective-augmentation-centralizer} gives the following result.

\begin{cor}\label{thm:aug-conj-qdl}
Let $S\subset G$ be a set of pairwise non-conjugate elements. Suppose that
$C_G(s)=\langle s\rangle$ for every $s\in S$. Then the augmented quandle
\[
Q=\bigcup_{g\in G}g^{-1}Sg\hookrightarrow G
\]
is free on the inclusion $S\hookrightarrow Q$.
\end{cor}

\begin{proof}
The inclusion $Q\hookrightarrow G$ is injective. Hence Corollary~\ref{cor:injective-augmentation-centralizer} gives
$\Stab_G(s)=C_G(s)=\langle s\rangle$. Proposition~\ref{pro:free-aq-coset-criterion} now gives the result.
\end{proof}

\begin{eg}\label{eg:hyperbolic-conjugation-quandle}
Let $G$ be a torsion-free hyperbolic group and let $g\in G$ be a nontrivial element. 
Then $\langle g\rangle$ has finite index in $C_G(g)$.
See Bridson--Haefliger~\cite{MR1744486}*{Corollary~III.$\Gamma$.3.10}. 
Thus $C_G(g)$ is virtually cyclic. Since $C_G(g)$ is torsion-free, it is infinite cyclic. If $g$ is not a proper power, then
$C_G(g)=\langle g\rangle$.

Let $S\subset G$ be a set of pairwise non-conjugate elements, each of which is nontrivial and not a proper power. Then Corollary~\ref{thm:aug-conj-qdl} shows that
\[
\partial\colon\bigcup_{g\in G}g^{-1}Sg\hookrightarrow G
\]
is a free augmented quandle with basis $S$.
The associated crossed module is the free crossed module on the inclusion $S\hookrightarrow G$.
\end{eg}

\begin{eg}\label{eg:symm_gp_n-cycle}
Let $n\geq3$, and let $\mathfrak S_n$ and $\mathfrak A_n$ denote the symmetric and alternating groups on $n$ letters, respectively.
The inclusion of the conjugacy class of an $n$-cycle in
$\mathfrak S_n$ is a free augmented quandle with a singleton basis. Indeed,
if $\sigma$ is an $n$-cycle, then
\[
C_{\mathfrak S_n}(\sigma)=\langle\sigma\rangle.
\]
The same conclusion holds for the conjugacy class of an $(n-1)$-cycle.
The associated crossed module is the free crossed module on the map
\[
\{*\}\to\mathfrak S_n
\]
sending $*$ to the chosen $n$-cycle or $(n-1)$-cycle. Its image is the normal closure of that cycle in $\mathfrak S_n$, and hence is $\mathfrak S_n$ when the cycle is odd and $\mathfrak A_n$ when it is even.
\end{eg}

\subsection{Generalized Alexander quandles}
Let $M=\Alex(G,\varphi)$, where $G$ is a group and
$\varphi\in\Aut(G)$. Thus the underlying set of $M$ is $G$, with operation
\[
g\ast h=\varphi(gh^{-1})h.
\]

\begin{pro}[Bonatto~\cite{Bonatto2020}*{Lemma~1.7 and Section~2.1}]
\label{pro:Inn-generalized-Alexander}
Put
\[
\Gamma\coloneqq
\langle \varphi(g)^{-1}g\mid g\in G\rangle.
\]
Then $\varphi(\Gamma)=\Gamma$. Let $\langle\varphi\rangle$ act on $\Gamma$ from the right by
$a\cdot\varphi^n=\varphi^n(a)$.
Then
\[
\Inn(M)\cong \langle\varphi\rangle\ltimes\Gamma,
\]
where
\[
(\varphi^n,a)(\varphi^m,b)
=
(\varphi^{n+m},\varphi^m(a)b).
\]
Under this identification, $(\varphi^n,a)$ acts on $M=G$ by
\[
T_{a,\varphi^n}(g)=\varphi^n(g)a.
\]
\end{pro}

We use this description of $\Inn(M)$ to verify the freeness criterion for the canonical inner augmentation.
\begin{thm}\label{thm:generalized-Alexander-free}
Let $M=\Alex(G,\varphi)$, and let $S\subset M$ be a complete set of representatives for the $\Inn(M)$-orbits. Then the canonical inner augmentation
\[
M\to\Inn(M),
\quad
x\mapsto S_x
\]
is a free augmented quandle with basis $S$.
\end{thm}

\begin{proof}
We apply Proposition~\ref{pro:free-aq-coset-criterion}. By Proposition~\ref{pro:Inn-generalized-Alexander}, every element of $\Inn(M)$ has the form
\[
T_{a,\varphi^n}(g)=\varphi^n(g)a
\quad(a\in\Gamma,\ n\in\mathbb Z).
\]
Fix $x\in M$. The element $T_{a,\varphi^n}$ fixes $x$ if and only if
\[
a=\varphi^n(x)^{-1}x.
\]
On the other hand, for every $n\in\mathbb Z$,
\[
S_x^n(g)=\varphi^n(g)\varphi^n(x)^{-1}x.
\]
Thus the elements of $\Inn(M)$ fixing $x$ are precisely the powers of $S_x$. Hence
\[
\Stab_{\Inn(M)}(x)=\langle S_x\rangle.
\]
The result follows from Proposition~\ref{pro:free-aq-coset-criterion}.
\end{proof}

\subsection{Knot quandles}
We finish with applications to knot quandles. For knots $K\subset S^3$, we use Joyce's coset description~\cite{MR0638121}. More generally, for an $n$-knot $K\subset S^{n+2}$, the isomorphism
$\As(Q_K)\cong G_K$ follows from Fenn--Rourke~\cite{MR1194995}*{Propositions~3.1 and~3.2}. See also Winter~\cite{MR3434551}*{Theorem~10} for the reconstruction of the knot group and peripheral structure from the knot quandle in arbitrary dimensions.

Let $P_K\le G_K$ denote the peripheral subgroup. In the standard coset description
\[
Q_K\cong P_K\backslash G_K,
\]
the stabilizer of the coset $P_K$ under the $G_K$-action is precisely $P_K$. Under the canonical homomorphism
\[
G_K\to\Inn(Q_K),
\]
the subgroup $P_K$ maps onto the stabilizer of the coset $P_K\in P_K\backslash G_K\cong Q_K$. For the classical case, see Joyce~\cite{MR0638121}. For the description in terms of the peripheral subgroup in arbitrary dimensions, see Winter~\cite{MR3434551}*{Lemma~1 and Theorem~8}.

Fix a meridian $m\in P_K$. 
Under the identification above, the standard knot augmentation is given by
\[
\varepsilon_K\colon Q_K\to G_K,
\quad
P_Kg\mapsto g^{-1}mg.
\]
The $G_K$-action on $Q_K$ is transitive. Since this action factors through the canonical homomorphism
$G_K\to\Inn(Q_K)$, the quandle $Q_K$ is connected.

\begin{eg}\label{eg:knot-converse-failure}
Let $K\subset S^3$ be a nontrivial knot, let $Q_K$ be its knot quandle, and let $G_K$ be the knot group.
Choose any $q\in Q_K$ and put
$\mu=\varepsilon_K(q)$. Then $\mu$ is a meridian. Under the classical isomorphism
$\As(Q_K)\cong G_K$, the associated crossed module
$\As(Q_K)\to G_K$ is identified with
$\mathrm{id}\colon G_K\to G_K$. By Example~\ref{eg:knot-group-identity-free}, this associated crossed module is free with basis $\{e_q\}$.

On the other hand, the stabilizer of $q$ under the $G_K$-action is conjugate to the peripheral subgroup $P_K\cong\mathbb Z^2$, whereas $\langle\mu\rangle$ is cyclic. Proposition~\ref{pro:free-aq-coset-criterion} therefore shows that
$\varepsilon_K\colon Q_K\to G_K$ is not a free augmented quandle. 
Thus every nontrivial classical knot gives an augmented quandle that is not free although its associated crossed module is free. 
\end{eg}

\begin{thm}\label{thm:one-knot-free-augmented-quandle}
Let $K\subset S^3$ be a knot and let $Q_K$ be its knot quandle. The canonical inner augmentation
\[
Q_K\to\Inn(Q_K)
\]
is a free augmented quandle with a singleton basis if and only if $K$ is a torus knot. Here the unknot is regarded as a trivial torus knot. Since $Q_K$ is connected, when the augmentation is free, any $q\in Q_K$ may be chosen as a basis element.
\end{thm}

\begin{proof}
For the unknot, $Q_K$ is the one-element quandle, so the statement is immediate. Assume that $K$ is nontrivial.
By Joyce's description~\cite{MR0638121}, the elements of $G_K$ acting trivially on $Q_K$ under the standard right action associated with the augmentation $Q_K\to G_K$ are precisely the elements of $Z(G_K)$.
Hence
\[
\Inn(Q_K)\cong G_K/Z(G_K).
\]
Choose $x\in Q_K$ corresponding to the coset $P_K$ under the identification
$Q_K\cong P_K\backslash G_K$, and put $\mu=\varepsilon_K(x)$. Then $\mu$ is a meridian, and the stabilizer of $x$ in $\Inn(Q_K)$ is the image of $P_K$ in $G_K/Z(G_K)$.

Suppose that $K=T(p,r)$ is a nontrivial torus knot. Write
\[
P_K=\langle\mu,\lambda\rangle,
\]
where $\lambda$ is the preferred longitude. There is a generator $z$ of $Z(G_K)$ such that
\[
\lambda=z\mu^{-pr}.
\]
Hence the image of $P_K$ in $G_K/Z(G_K)$ is generated by the image of $\mu$, which corresponds to $S_x\in\Inn(Q_K)$. Thus
\[
\Stab_{\Inn(Q_K)}(x)=\langle S_x\rangle.
\]
Proposition~\ref{pro:free-aq-coset-criterion} shows that the canonical inner augmentation is free. Since $\Inn(Q_K)$ acts transitively on $Q_K$, the same stabilizer equality holds for every $q\in Q_K$, so any $q\in Q_K$ may be chosen as a basis element.

If $K$ is not a torus knot, then $Z(G_K)=1$. Hence, for every $q\in Q_K$, the stabilizer of $q$ in $\Inn(Q_K)$ is conjugate to $P_K\cong\mathbb Z^2$, and therefore cannot equal the cyclic group $\langle S_q\rangle$. Proposition~\ref{pro:free-aq-coset-criterion} shows that the canonical inner augmentation is not free. For these classical facts, see Burde--Zieschang--Heusener~\cite{BurdeZieschangHeusener2014}.
\end{proof}

\begin{rem}\label{rem:torus-versus-nontorus}
Theorem~\ref{thm:one-knot-free-augmented-quandle} highlights the following contrast.
If $K$ is a non-torus knot, then $Z(G_K)=1$ and hence
$\Inn(Q_K)\cong G_K$, but the augmentation $Q_K\to\Inn(Q_K)$ is not free. If $K=T(p,q)$ is a nontrivial torus knot, then the augmentation is free, whereas
\[
\Inn(Q_K)\cong G_K/Z(G_K)\not\cong G_K.
\]
Indeed,
\[
G_K^{\mathrm{ab}}\cong\mathbb Z,
\qquad
(G_K/Z(G_K))^{\mathrm{ab}}
\cong C_{|p|}\oplus C_{|q|},
\]
where $C_n$ denotes the cyclic group of order $n$, so the two groups are not isomorphic.
\end{rem}

\begin{cor}\label{cor:torus-knot-free-crossed-module}
Let $K=T(p,q)$ be a nontrivial torus knot with coprime integers $p$ and $q$, $|p|,|q|>1$, and let
\[
G_K=\langle a,b\mid a^p=b^q\rangle
\]
be its knot group. Then the projection
\[
\partial\colon G_K\to\Inn(Q_K)\cong G_K/Z(G_K)
\]
is a free crossed module with basis $\{\mu\}$, where $\mu$ is a meridian. Equivalently, since
$Z(G_K)=\langle a^p\rangle=\langle b^q\rangle$, the projection
\[
\partial\colon \langle a,b\mid a^p=b^q\rangle
\to
\langle a,b\mid a^p=b^q=1\rangle
\]
is a free crossed module. If $s,r\in\mathbb Z$ satisfy $ps-qr=1$, one may take the meridian to be $a^{-r}b^s$.
\end{cor}

\begin{proof}
Choose $q_0\in Q_K$ whose image in $G_K$ is a meridian $\mu$. By Theorem~\ref{thm:one-knot-free-augmented-quandle}, the augmentation
$Q_K\to\Inn(Q_K)$ is a free augmented quandle with basis $\{q_0\}$. By Theorem~\ref{thm:AT3}, its associated crossed module is free with basis $\{e_{q_0}\}$. Under Joyce's classical identification
$\As(Q_K)\cong G_K$~\cite{MR0638121}, the element $e_{q_0}$ corresponds to $\mu$. Under the identifications
\[
\As(Q_K)\cong G_K,
\qquad
\Inn(Q_K)\cong G_K/Z(G_K),
\]
the boundary map of the associated crossed module is the canonical projection
$G_K\to G_K/Z(G_K)$. 
The presentation of $G_K/Z(G_K)$ follows from the standard presentation of torus knot groups and the description of their center.
See Burde--Zieschang--Heusener~\cite{BurdeZieschangHeusener2014}.
\end{proof}

\begin{rem}
Corollary~\ref{cor:torus-knot-free-crossed-module} can also be verified directly from Ratcliffe's criterion. Put
\[
C=G_K,
\qquad
G=G_K/Z(G_K)\cong C_{|p|}*C_{|q|},
\qquad
N=G.
\]
Then $G/N=1$ and $C_{\mathrm{ab}}\cong\mathbb Z$ is free with basis the image of a meridian. The image of a meridian normally generates $G$, and
\[
H_2(G)=H_2(C_{|p|}*C_{|q|})=0.
\]
Thus the three conditions in Theorem~\ref{thm:Ratcliffe} are satisfied.
\end{rem}

\begin{thm}\label{thm:higher-knot-free-augmented-quandle}
Let $n\geq2$ and let $K\subset S^{n+2}$ be an $n$-knot. Then the canonical inner augmentation
\[
Q_K\to\Inn(Q_K)
\]
is a free augmented quandle with a singleton basis. More precisely, $Q_K$ is connected, and any $q\in Q_K$ may be chosen as a basis element. Moreover, the associated crossed module
\[
\As(Q_K)\to\Inn(Q_K)
\]
is free with basis $\{e_q\}$. Under the isomorphism $\As(Q_K)\cong G_K$, the element $e_q$ corresponds to a conjugate of a meridian.
\end{thm}

\begin{proof}
Choose any $q\in Q_K$ and write $q=P_Kg$ under the coset description. Put $\mu=g^{-1}mg$, where $m$ is a meridian. By the description of the knot quandle in terms of the peripheral subgroup, the stabilizer of $q$ in $\Inn(Q_K)$ is the image of $g^{-1}P_Kg$.
See Fenn--Rourke~\cite{MR1194995} and Winter~\cite{MR3434551}*{Lemma~1 and Theorem~8}.
For $n\geq2$, since $\partial N(K)\cong S^n\times S^1$ for a tubular neighborhood $N(K)$ of $K$, we have
\[
P_K=\langle m\rangle,
\]
and hence
\[
g^{-1}P_Kg=\langle\mu\rangle.
\]
Since $\mu=\varepsilon_K(q)$ maps to $S_q$ under $G_K\to\Inn(Q_K)$, it follows that
\[
\Stab_{\Inn(Q_K)}(q)=\langle S_q\rangle.
\]
Proposition~\ref{pro:free-aq-coset-criterion} shows that the canonical inner augmentation is free with basis $\{q\}$. The statement about the associated crossed module follows from Theorem~\ref{thm:AT3}.
\end{proof}

\begin{rem}\label{rem:knot-peripheral-dimension}
The classical and higher-dimensional cases differ in their peripheral subgroups. Let $P_K\le G_K$ be the peripheral subgroup and let $m$ be a meridian. 
For a classical knot,
\[
P_K=\langle m,\lambda\rangle,
\]
where $\lambda$ is the preferred longitude.
If the knot is nontrivial, then $P_K\cong\mathbb Z^2$. For an $n$-knot with $n\geq2$,
\[
P_K=\langle m\rangle.
\]
Thus, for $n\geq2$, the stabilizer condition in Proposition~\ref{pro:free-aq-coset-criterion} is automatic. In the classical case, the passage from $G_K$ to
$\Inn(Q_K)\cong G_K/Z(G_K)$ is essential. With the convention that the unknot is a trivial torus knot,
\[
Z(G_K)\ne1
\quad\Longleftrightarrow\quad
K\text{ is a torus knot}.
\]
For a nontrivial torus knot, the image of $P_K$ in $G_K/Z(G_K)$ is the cyclic subgroup generated by the image of the meridian. For a non-torus knot, the center is trivial and the stabilizer remains $\mathbb Z^2$.
\end{rem}

\section*{Acknowledgement}
The first author was partially supported by JSPS KAKENHI Grant Number 24K06727.

\section*{Tool Use}
OpenAI's ChatGPT was used for English-language editing,
LaTeX assistance, reference checking, and suggestions
concerning the presentation and consistency of several
arguments.
All mathematical content was reviewed and verified by the authors.


\begin{bibdiv}
\begin{biblist}

\bib{MR4175808}{article}{
   author={Akita, Toshiyuki},
   title={The adjoint group of a Coxeter quandle},
   journal={Kyoto J. Math.},
   volume={60},
   date={2020},
   number={4},
   pages={1245--1260},
   issn={2156-2261},
   review={\MR{4175808}},
   doi={10.1215/21562261-2019-0061},
}

\bib{MR4564617}{article}{
   author={Akita, Toshiyuki},
   title={Embedding Alexander quandles into groups},
   journal={J. Knot Theory Ramifications},
   volume={32},
   date={2023},
   number={2},
   pages={Paper No. 2350011, 4},
   issn={0218-2165},
   review={\MR{4564617}},
   doi={10.1142/S0218216523500116},
}

\bib{AkitaTakaseOCAMI2025}{article}{
   author={Akita, Toshiyuki},
   author={Takase, Sota},
   title={Freeness of crossed modules, augmented quandles, and augmented racks},
   book={
      title={Quandles and Symmetric Spaces 2025},
      series={OCAMI Reports},
      volume={7},
      publisher={Osaka Central Advanced Mathematical Institute (OCAMI),
         Osaka Metropolitan University},
      date={2025},
   },
   pages={113--127},
}

\bib{MR2799090}{article}{
   author={Andruskiewitsch, N.},
   author={Fantino, F.},
   author={Garc{\'{\i}}a, G. A.},
   author={Vendramin, L.},
   title={On Nichols algebras associated to simple racks},
   conference={
      title={Groups, algebras and applications},
   },
   book={
      series={Contemp. Math.},
      volume={537},
      publisher={Amer. Math. Soc., Providence, RI},
   },
   date={2011},
   pages={31--56},
   review={\MR{2799090 (2012g:16065)}},
   doi={10.1090/conm/537/10565},
}

\bib{MR1994219}{article}{
   author={Andruskiewitsch, Nicol\'as},
   author={Gra\~na, Mat\'ias},
   title={From racks to pointed Hopf algebras},
   journal={Adv. Math.},
   volume={178},
   date={2003},
   number={2},
   pages={177--243},
   issn={0001-8708},
   review={\MR{1994219}},
   doi={10.1016/S0001-8708(02)00071-3},
}


\bib{MR2825807}{article}{
   author={Baez, John C.},
   author={Huerta, John},
   title={An invitation to higher gauge theory},
   journal={Gen. Relativity Gravitation},
   volume={43},
   date={2011},
   number={9},
   pages={2335--2392},
   issn={0001-7701},
   review={\MR{2825807}},
   doi={10.1007/s10714-010-1070-9},
}


\bib{MR2342821}{article}{
   author={Baez, John C.},
   author={Schreiber, Urs},
   title={Higher gauge theory},
   conference={
      title={Categories in algebra, geometry and mathematical physics},
   },
   book={
      series={Contemp. Math.},
      volume={431},
      publisher={Amer. Math. Soc., Providence, RI},
   },
   isbn={978-0-8218-3970-6},
   isbn={0-8218-3970-5},
   date={2007},
   pages={7--30},
   review={\MR{2342821}},
   doi={10.1090/conm/431/08264},
}


\bib{arXiv:2506.23175}{book}{
   author={Bardakov, Valeriy G.},
   author={Elhamdadi, Mohamed},
   author={Singh, Mahender},
   title={Algebra of the Yang--Baxter Equation:
          Skew Braces, Quandles, and Cohomology},
   series={Springer Monographs in Mathematics},
   publisher={Springer Nature Singapore},
   date={2026},
   isbn={978-981-95-8929-6},
}

\bib{BardakovNasybullov2019}{article}{
   author={Bardakov, Valeriy},
   author={Nasybullov, Timur},
   title={Embeddings of quandles into groups},
   journal={J. Algebra Appl.},
   volume={19},
   date={2020},
   number={7},
   pages={2050136, 20},
   issn={0219-4988},
   review={\MR{4129183}},
   doi={10.1142/S0219498820501364},
}

\bib{Bonatto2020}{article}{
   author={Bonatto, Marco},
   title={Principal and doubly homogeneous quandles},
   journal={Monatsh. Math.},
   volume={191},
   date={2020},
   number={4},
   pages={691--717},
   issn={0026-9255},
   review={\MR{4075376}},
   doi={10.1007/s00605-019-01334-1},
}



\bib{MR1744486}{book}{
   author={Bridson, Martin R.},
   author={Haefliger, Andr\'e},
   title={Metric spaces of non-positive curvature},
   series={Grundlehren der mathematischen Wissenschaften [Fundamental
   Principles of Mathematical Sciences]},
   volume={319},
   publisher={Springer-Verlag, Berlin},
   date={1999},
   pages={xxi+643},
   isbn={3-540-64324-9},
   review={\MR{1744486}},
   doi={10.1007/978-3-662-12494-9},
}

\bib{MR975077}{article}{
   author={Brieskorn, E.},
   title={Automorphic sets and braids and singularities},
   conference={
      title={Braids},
      address={Santa Cruz, CA},
      date={1986},
   },
   book={
      series={Contemp. Math.},
      volume={78},
      publisher={Amer. Math. Soc., Providence, RI},
   },
   date={1988},
   pages={45--115},
   review={\MR{975077}},
   doi={10.1090/conm/078/975077},
}

\bib{MR0579818}{article}{
   author={Brown, Ronald},
   title={On the second relative homotopy group of an adjunction space: an
   exposition of a theorem of J. H. C. Whitehead},
   journal={J. London Math. Soc. (2)},
   volume={22},
   date={1980},
   number={1},
   pages={146--152},
   issn={0024-6107},
   review={\MR{0579818}},
   doi={10.1112/jlms/s2-22.1.146},
}

\bib{MR0478150}{article}{
   author={Brown, Ronald},
   author={Higgins, Philip J.},
   title={On the connection between the second relative homotopy groups of
   some related spaces},
   journal={Proc. London Math. Soc. (3)},
   volume={36},
   date={1978},
   number={2},
   pages={193--212},
   issn={0024-6115},
   review={\MR{0478150}},
   doi={10.1112/plms/s3-36.2.193},
}

\bib{MR2841564}{book}{
   author={Brown, Ronald},
   author={Higgins, Philip J.},
   author={Sivera, Rafael},
   title={Nonabelian algebraic topology},
   series={EMS Tracts in Mathematics},
   volume={15},
   note={Filtered spaces, crossed complexes, cubical homotopy groupoids;
   With contributions by Christopher D. Wensley and Sergei V. Soloviev},
   publisher={European Mathematical Society (EMS), Z\"urich},
   date={2011},
   pages={xxxvi+668},
   isbn={978-3-03719-083-8},
   review={\MR{2841564}},
   doi={10.4171/083},
}

\bib{BrownHuebschmann1982}{article}{
   author={Brown, Ronald},
   author={Huebschmann, Johannes},
   title={Identities among relations},
   conference={
      title={Low-dimensional topology},
      address={Bangor},
      date={1979},
   },
   book={
      series={London Math. Soc. Lecture Note Ser.},
      volume={48},
      publisher={Cambridge Univ. Press, Cambridge},
   },
   date={1982},
   pages={153--202},
   doi={10.1017/CBO9780511758935.010},
}

\bib{MR0419643}{article}{
   author={Brown, Ronald},
   author={Spencer, Christopher B.},
   title={$G$-groupoids, crossed modules and the fundamental groupoid of a
   topological group},
   note={Nederl. Akad. Wetensch. Proc. Ser. A {\bf 79}},
   journal={Indag. Math.},
   volume={38},
   date={1976},
   number={4},
   pages={296--302},
   review={\MR{0419643}},
   doi={10.1016/1385-7258(76)90068-8},
}

\bib{BurdeZieschangHeusener2014}{book}{
   author={Burde, Gerhard},
   author={Zieschang, Heiner},
   author={Heusener, Michael},
   title={Knots},
   series={De Gruyter Studies in Mathematics},
   volume={5},
   edition={extended edition},
   publisher={De Gruyter, Berlin},
   date={2014},
   pages={xiv+417},
   isbn={978-3-11-027074-7},
   isbn={978-3-11-027078-5},
   review={\MR{3156509}},
}

\bib{CransWagemann2014}{article}{
   author={Crans, Alissa S.},
   author={Wagemann, Friedrich},
   title={Crossed modules of racks},
   journal={Homology Homotopy Appl.},
   volume={16},
   date={2014},
   number={2},
   pages={85--106},
   doi={10.4310/HHA.2014.v16.n2.a5},
}

\bib{MR1279178}{article}{
   author={Dyer, Micheal N.},
   title={Crossed modules and $\Pi_2$ homotopy modules},
   conference={
      title={Two-dimensional homotopy and combinatorial group theory},
   },
   book={
      series={London Math. Soc. Lecture Note Ser.},
      volume={197},
      publisher={Cambridge Univ. Press, Cambridge},
   },
   isbn={0-521-44700-3},
   date={1993},
   pages={125--156},
   review={\MR{1279178}},
   doi={10.1017/CBO9780511629358.006},
}

\bib{MR3971587}{book}{
   author={Ellis, Graham},
   title={An invitation to computational homotopy},
   publisher={Oxford University Press, Oxford},
   date={2019},
   pages={xx+525},
   isbn={978-0-19-883298-0},
   isbn={978-0-19-883297-3},
   review={\MR{3971587}},
   doi={10.1093/oso/9780198832973.001.0001},
}


\bib{MR0825179}{article}{
   author={Ellis, G. J.},
   author={Porter, T.},
   title={Free and projective crossed modules and the second homology group
   of a group},
   journal={J. Pure Appl. Algebra},
   volume={40},
   date={1986},
   number={1},
   pages={27--31},
   issn={0022-4049},
   review={\MR{0825179}},
   doi={10.1016/0022-4049(86)90028-9},
}


\bib{MR1194995}{article}{
   author={Fenn, Roger},
   author={Rourke, Colin},
   title={Racks and links in codimension two},
   journal={J. Knot Theory Ramifications},
   volume={1},
   date={1992},
   number={4},
   pages={343--406},
   issn={0218-2165},
   review={\MR{1194995}},
   doi={10.1142/S0218216592000203},
}

\bib{MR1364012}{article}{
   author={Fenn, Roger},
   author={Rourke, Colin},
   author={Sanderson, Brian},
   title={Trunks and classifying spaces},
   journal={Appl. Categ. Structures},
   volume={3},
   date={1995},
   number={4},
   pages={321--356},
   issn={0927-2852},
   review={\MR{1364012}},
   doi={10.1007/BF00872903},
}

\bib{MR1020583}{article}{
   author={Freyd, Peter J.},
   author={Yetter, David N.},
   title={Braided compact closed categories with applications to
   low-dimensional topology},
   journal={Adv. Math.},
   volume={77},
   date={1989},
   number={2},
   pages={156--182},
   issn={0001-8708},
   review={\MR{1020583}},
   doi={10.1016/0001-8708(89)90018-2},
}

\bib{MR0816893}{article}{
   author={Guti\'errez, Mauricio A.},
   author={Hirschhorn, Philip S.},
   title={Free simplicial groups and the second relative homotopy group of
   an adjunction space},
   journal={J. Pure Appl. Algebra},
   volume={39},
   date={1986},
   number={1-2},
   pages={119--123},
   issn={0022-4049},
   review={\MR{0816893}},
   doi={10.1016/0022-4049(86)90139-8},
}

\bib{MR4322147}{article}{
   author={Huebschmann, Johannes},
   title={On the history of Lie brackets, crossed modules, and Lie-Rinehart
   algebras},
   journal={J. Geom. Mech.},
   volume={13},
   date={2021},
   number={3},
   pages={385--402},
   issn={1941-4889},
   review={\MR{4322147}},
   doi={10.3934/jgm.2021009},
}

\bib{MR4666020}{article}{
   author={Huebschmann, Johannes},
   title={Crossed modules},
   journal={Notices Amer. Math. Soc.},
   volume={70},
   date={2023},
   number={11},
   pages={1802--1813},
   issn={0002-9920},
   review={\MR{4666020}},
   doi={10.1090/noti2831},
}

\bib{MR2628474}{book}{
   author={Joyce, David},
   title={An algebraic approach to symmetry with applications to knot theory},
   note={Thesis (Ph.D.)--University of Pennsylvania},
   publisher={ProQuest LLC, Ann Arbor, MI},
   date={1979},
   pages={127},
   review={\MR{2628474}},
}

\bib{MR0638121}{article}{
   author={Joyce, David},
   title={A classifying invariant of knots, the knot quandle},
   journal={J. Pure Appl. Algebra},
   volume={23},
   date={1982},
   number={1},
   pages={37--65},
   issn={0022-4049},
   review={\MR{0638121}},
   doi={10.1016/0022-4049(82)90077-9},
}

\bib{Kohno2020}{article}{
   author={Kohno, Toshitake},
   title={Higher holonomy maps for hyperplane arrangements},
   journal={Eur. J. Math.},
   volume={6},
   date={2020},
   number={3},
   pages={905--927},
   issn={2199-675X},
   review={\MR{4151724}},
   doi={10.1007/s40879-019-00382-z},
}

\bib{MR0672410}{article}{
   author={Matveev, S. V.},
   title={Distributive groupoids in knot theory},
   language={Russian},
   journal={Mat. Sb. (N.S.)},
   volume={119(161)},
   date={1982},
   number={1},
   pages={78--88, 160},
   issn={0368-8666},
   review={\MR{0672410}},
}

\bib{nosaka-book}{book}{
   author={Nosaka, Takefumi},
   title={Quandles and topological pairs},
   series={SpringerBriefs in Mathematics},
   note={Symmetry, knots, and cohomology},
   publisher={Springer, Singapore},
   date={2017},
   pages={ix+136},
   isbn={978-981-10-6792-1},
   isbn={978-981-10-6793-8},
   review={\MR{3729413}},
   doi={10.1007/978-981-10-6793-8},
}

\bib{arXiv:2610.38066}{article}{
   author={Nosaka, Takefumi},
   title={Loop space models and integral homology of connected quandles},
   date={2026},
   eprint={https://arxiv.org/abs/2610.38066},
}

\bib{MR1808345}{article}{
   author={Oda, Fumihito},
   author={Yoshida, Tomoyuki},
   title={Crossed Burnside rings. I. The fundamental theorem},
   journal={J. Algebra},
   volume={236},
   date={2001},
   number={1},
   pages={29--79},
   issn={0021-8693},
   review={\MR{1808345}},
   doi={10.1006/jabr.2000.8341},
}

\bib{MR2095572}{article}{
   author={Oda, Fumihito},
   author={Yoshida, Tomoyuki},
   title={Crossed Burnside rings. II. The Dress construction of a Green
   functor},
   journal={J. Algebra},
   volume={282},
   date={2004},
   number={1},
   pages={58--82},
   issn={0021-8693},
   review={\MR{2095572}},
   doi={10.1016/j.jalgebra.2004.07.040},
}

\bib{MR2746028}{article}{
   author={Oda, Fumihito},
   author={Yoshida, Tomoyuki},
   title={The crossed Burnside rings III: the dress construction for a
   Tambara functor},
   journal={J. Algebra},
   volume={327},
   date={2011},
   pages={31--49},
   issn={0021-8693},
   review={\MR{2746028}},
   doi={10.1016/j.jalgebra.2010.10.018},
}

\bib{MR0579810}{article}{
   author={Ratcliffe, John G.},
   title={Free and projective crossed modules},
   journal={J. London Math. Soc. (2)},
   volume={22},
   date={1980},
   number={1},
   pages={66--74},
   issn={0024-6107},
   review={\MR{0579810}},
   doi={10.1112/jlms/s2-22.1.66},
}

\bib{Takase2024}{thesis}{
   author={Takase, Sota},
   title={Free objects in the category of augmented quandles and groups admitting Wirtinger presentations (in Japanese)},
   type={Master's thesis},
   organization={Hokkaido University},
   date={2024},
}

\bib{MR4368859}{book}{
   author={Wagemann, Friedrich},
   title={Crossed modules},
   series={De Gruyter Studies in Mathematics},
   volume={82},
   publisher={De Gruyter, Berlin},
   date={2021},
   pages={xiv+393},
   isbn={978-3-11-075076-8},
   isbn={978-3-11-075095-9},
   isbn={978-3-11-075099-7},
   review={\MR{4368859}},
   doi={10.1515/9783110750959},
}

\bib{MR0004123}{article}{
   author={Whitehead, J. H. C.},
   title={On adding relations to homotopy groups},
   journal={Ann. of Math. (2)},
   volume={42},
   date={1941},
   pages={409--428},
   issn={0003-486X},
   review={\MR{0004123}},
   doi={10.2307/1968907},
}

\bib{MR0017537}{article}{
   author={Whitehead, J. H. C.},
   title={Note on a previous paper entitled ``On adding relations to
   homotopy groups.''},
   journal={Ann. of Math. (2)},
   volume={47},
   date={1946},
   pages={806--810},
   issn={0003-486X},
   review={\MR{0017537}},
   doi={10.2307/1969237},
}

\bib{MR0030760}{article}{
   author={Whitehead, J. H. C.},
   title={Combinatorial homotopy. II},
   journal={Bull. Amer. Math. Soc.},
   volume={55},
   date={1949},
   pages={453--496},
   issn={0002-9904},
   review={\MR{0030760}},
   doi={10.1090/S0002-9904-1949-09213-3},
}

\bib{MR3434551}{article}{
   author={Winter, Blake K.},
   title={Virtual links in arbitrary dimensions},
   journal={J. Knot Theory Ramifications},
   volume={24},
   date={2015},
   number={14},
   pages={1550062, 38},
   issn={0218-2165},
   review={\MR{3434551}},
   doi={10.1142/S0218216515500625},
}

\end{biblist}
\end{bibdiv}
\end{document}